\documentclass[12 pt]{article}

\usepackage{amsmath,amssymb,amsthm}
\usepackage[T1]{fontenc}
\usepackage{amsfonts}
\usepackage{epsfig}
\usepackage{graphicx}
\usepackage{epic}
\usepackage{authblk}
\usepackage{color}
\usepackage{psfrag}
\usepackage{verbatim}
\usepackage{tikz}
\usepackage{psfrag}
\usepackage{caption}
\usepackage{subcaption}
\providecommand{\keywords}[1]{\textbf{\textit{Keywords---}} #1}

\title{Flux-conservative Hermite methods for simulation of nonlinear conservation laws \thanks{Supported in part by NSF grant DMS-1319054. Any conclusions or recommendations expressed in this paper are those of the authors and do not necessarily reflect the views of NSF.}}

\author{Adeline Kornelus \thanks{kornelus@unm.edu} and Daniel Appel\"{o}\thanks{appelo@math.unm.edu}}  

\affil{Department of Mathematics and Statistics\\
  University of New Mexico, Albuquerque,NM}

\date{\today}

\begin{document}

\maketitle
\begin{abstract}
A new class of Hermite methods for solving nonlinear conservation laws is presented. While preserving the high order spatial accuracy for smooth solutions in the existing Hermite methods, the new  methods come with better stability properties. Artificial viscosity in the form of the entropy viscosity method is added to capture shocks.

\end{abstract}

\keywords{Hermite, conservative, high order method, the entropy viscosity method}

\section{Introduction}
\label{sec:intro}
Conservation laws model physical systems arising in traffic flows, aircraft design, weather forecast, and fluid dynamics. Numerical methods for conservation laws ideally conserve quantities like mass or energy, and accurately capture various physical components of the solutions, from small smooth scales to shock waves. The presence of both smooth waves and shock waves, for example in shock-turbulence interaction creates a challenge to the simulation of nonlinear conservation laws. 

Shock waves typically appear in solutions to nonlinear conservation laws, and are characterized by very thin regions where the solution changes rapidly. Approximation of shocks and small waves is challenging as the small and large scales need to be solved simultaneously. Historically, low order finite volume and finite difference methods equipped with flux/slope limiters have been used to handle shocks, see for example the textbooks \cite{leveque2002finite,LeVeque:1992ud}. The drawback with low order methods is that they cannot accurately propagate small scales over long distances and as a result, today the research focus has gravitated towards high order accurate methods with shock capturing capability.  

Among high order methods, the weighted essentially non-oscillatory (WENO) method, \cite{Shu1988439,shu1998essentially,Liu1994200}, has proven to be a method popular among practitioners. WENO methods are still relatively dissipative which may be a drawback for turbulent simulations, \cite{martin2006bandwidth}. Also, discontinuous Galerkin methods combined with shock capturing, \cite{dumbser2014posteriori,krivodonova2007limiters}, or selectively added artificial viscosity, \cite{persson2006sub,klockner2011viscous,zingan2013implementation}, have received significant interest. The latter approach traces back to the artificial viscosity method by Neumann and Richtmyer, \cite{NeuRich50} and the popular streamline diffusion method, \cite{BROOKS1982199,johnson1990convergence}, which computes the viscosity based on the residual of the PDE.

In this work, we advocate the combination of a high order method and selectively added artificial viscosity. Specifically, we show how the entropy viscosity by Guermond and Pasquetti, \cite{Guermond2008801}, can be implemented in our new flux-conservative Hermite methods.

First introduced by Goodrich, Hagstrom, and Lorenz in \cite{goodrich2006hermite} for hyperbolic initial-boundary value problems, Hermite methods use the solution and its first $m$ derivatives in each coordinate to construct an approximate solution to the PDE. The formulation by Goodrich et al. computes the flux at the cell centers, which for nonlinear problems leads to discontinuous fluxes at cell interfaces. This discontinuity results in lost of conservation. 

To address the lack of conservation, we develop a new class of Hermite methods, which share the basic features with the method in \cite{goodrich2006hermite}, such as interpolation and time evolution, but differs in the computation of the flux function. In the flux-conservative Hermite methods, the numerical flux is computed at cell edges and then interpolated to cell center for time evolution, hence by construction, the numerical flux is continuous at cell interface. Additionally, for nonlinear problems, it is more efficient to use one step methods than the Taylor series approach in \cite{goodrich2006hermite}, see \cite{hagstrom2007experiments,icosahom2014}. Here we will use the standard Runge-Kutta method to evolve in time.

The rest of the paper is organized as follows. In Section \ref{sec:conservation}, we derive conservation laws and discrete conservation. Then, in Section \ref{sec:original_hermite}, we describe Hermite methods as first introduced by Goodrich et al. \cite{goodrich2006hermite}, followed by the description of the flux-conservative Hermite methods in Section \ref{sec:conservative_hermite}. For shock capturing capability, we implement the entropy viscosity method, which is explained in Section \ref{sec:EV_intro}. In Section \ref{sec:numerical_exp}, we present the results of simulation on Euler's equations, see \cite{icosahom2016} for results on Burgers' equation.   

\section{Conservation Laws}
\label{sec:conservation}
A scalar conservation law in one space dimension takes the form 
\begin{equation}
  \label{eq:conservation_PDE}
  \frac{\partial}{\partial t} u(x,t) +  \frac{\partial}{\partial x} f(u(x,t)) = 0,
\end{equation}
where $u(x,t)$ is the state variable at location $x$ and time $t$ and $f(u)$ is the flux, or the rate of flow, of the state variable $u$.
 
The derivation of conservation laws comes from the observation that at any given time $t$, the rate of change of the total quantity of the state variable $u$ over some interval $[a,b]$ must be equal to the net flux $f(u)$ into the interval through the endpoints. Mathematically, this can be expressed as   
\begin{equation}
  \label{eq:conservation_integral}
  \frac{d}{dt} \int_a^b u \: dx = f(u(x(a),t))-f(u(x(b),t)). 
\end{equation}

When approximating the solution to scalar conservation laws given by equation (\ref{eq:conservation_PDE}), the PDE is typically discretized on a grid consisting of $N_x$ cells where $x_0=a$ and $x_{N_x}=b$. It is desirable that the numerical method satisfies discrete conservation. If the reconstructed solution $u_h^j$ and flux $f_h^j$ on any cell $\mathcal{I}_j = [x_{j-1},x_j]$ satisfy the condition $f_h^j|_{x=x_j} = f_{h}^{j+1}|_{x=x_j}$, $j=1,\ldots,N_x-1$, we immediately find
  \begin{align}
    \int_a^b \frac{\partial u_h^j}{\partial t}  \: dx & = \sum_{j=1}^{N_x} \int_{x_{j-1}}^{x_{j}} \frac{\partial u_h^j}{\partial t} \: dx \nonumber \\ 
                                              & = \sum_{j=1}^{N_x} \int_{x_{j-1}}^{x_{j}} \frac{\partial}{\partial x}(-{f_h^j}) \: dx \nonumber \\ 
                                              & = \sum_{j=1}^{N_x} \left({f_h^j}|_{x=x_{j-1}}-{f_h^j}|_{x=x_{j}} \right)\nonumber \\ 
                                              & = {f}_h^j(u(a))-{f}_{h}^j(u(b)).
  \end{align}

The property that $f_h^j|_{x=x_j} = f_{h}^{j+1}|_{x=x_j}$ does not hold for the original Hermite methods, and our goal here is to design a new Hermite method with this property. Before describing our new method, we briefly describe the original method. 

\section{Hermite Methods}
\label{sec:original_hermite}
A Hermite method of order $2m+1$ approximates the solution to a PDE by an element-wise polynomial that has continuous spatial derivatives up to order $m$ at the element's interfaces. In Hermite methods, the degrees of freedom are function and spatial derivative values, or equivalently the coefficients of the Taylor polynomial at the cell center of each element. The evolution of the degrees of freedom on each element can be performed locally.

\subsection{Hermite Method in One Dimension}
\label{sec:original_hermite_1D}
Consider again equation (\ref{eq:conservation_PDE}) on the domain $D = [x_L,x_R]$. Let $\mathcal{G}^p$ and $\mathcal{G}^d$ be the primal grid and the dual grid, defined as
\begin{equation}
\mathcal{G}^p = \left\{x_j = x_L+jh_x \right\}, \: j = 0,\ldots,N_x, 
\end{equation}

\begin{equation}
\mathcal{G}^d =\left\{x_{j+1/2} = x_L+\left(j+\frac{1}{2}\right)h_x\right\}, \: j = 0,\ldots,N_x-1,
\end{equation}
where $h_x = (x_R-x_L)/N_x$ is the distance between two adjacent nodes. Let $u^p$ and $u^d$ be the approximations to the solution on the primal and dual grids, respectively.

At the initial time $t_n= t_0+n\Delta t$, we assume that the approximation $u^p$ on the primal grid, the global piecewise polynomial  
\begin{equation}
u^p(x,t_n) = \sum_{k=0}^{m} c_{k}(t_n) \left(x-x_{j}\right)^k, x \in \mathcal{I}^d_j=[x_{j-1/2},x_{j+1/2}],
\end{equation}
is known. Starting from time $t=t_n$ on the primal grid $\mathcal{G}^p$, we evolve the solution one full time step by:

{\it Reconstruction by Hermite interpolation:}\label{it:reconstruction_1d} We construct $u^d$, the global Hermite interpolation polynomial of degree $(2m+1)$, on the dual grid. That is, 
   \begin{equation}
     \label{eq:hermite_interpolant_1d}
     u^d(x,t_n) = \sum_{k=0}^{2m+1} c_{k}(t_n) \left(x-x_{j-1/2}\right)^k, x \in \mathcal{I}^p_j=[x_{j-1},x_{j}], 
   \end{equation} 
   where the coefficients $c_k(t_n)$ are uniquely determined by the interpolation conditions
\begin{equation}
  \frac{\partial^k}{\partial x^k}u^d(x_i,t_n) = \frac{d^k}{dx^k}u^p(x_i,t_n), \ \ k=0,\ldots,m,\ \ i=j-1,j.
\end{equation}

{\it Time evolution:}\label{it:time_evolution_1d} By rewriting equation (\ref{eq:conservation_PDE}) as $u_t = -f(u)_x$, it is obvious that in order to evolve $u^d$, we need to compute a polynomial approximation $f^d$ to the flux function $f(u)$. We offer two ways to obtain $f^d$:
  \begin{itemize}
  \item Modal approach: Perform polynomial operations (addition, substraction, multiplication, and division) on $u^d$ and truncate the resulting polynomial to suitable degree.
  \item Pseudospectral approach: Compute a local polynomial $f_h^*$ interpolating $f(u^d)$ on a quadrature grid $\mathcal{G}^{ps}$ inside $\mathcal{I}^p_j, j=1,\ldots,N_x$, and transform $f_h^*$ to a Taylor polynomial $f^d$. 
  \end{itemize}
  We differentiate $f^d$ in polynomial sense to get an approximation to the derivative of the flux function $f(u)_x$. We usually use the modal approach unless this option is not applicable. Next, we derive an ODE to evolve $u^d$, or equivalently the coefficients of the Hermite interpolant ${\mathbf{c}(t)=(c_0(t),\ldots,c_{2m+1}(t))^T}$, by insisting that the numerical solution $u^d$ satisfy equation (\ref{eq:conservation_PDE}) along with derivatives of (\ref{eq:conservation_PDE}), at the cell centers $x=x_{j-1/2}, j=1,\ldots,N_x$. The resulting system of ODE for $c_k, k=0,\ldots,2m+1$, can then be evolved independently on each $\mathcal{I}^p_j$ with any ODE solver. The reconstruction step provides the initial data, $\mathbf{c}(t_n)$. 
     
{\it Repeat on dual grid:} At the end of the half time step, we have evolved the degree $2m+1$ polynomial $u^d$ from time $t=t_n$ to $t=t_{n+1/2}$. Before repeating the above process, we truncate $u^d(x,t_{n+1/2})$ by removing the coefficients $c_k$, $k=m+1,\ldots,2m+1$. We then repeat the above process (including the truncation) to obtain $u^p$ at time $t=t_{n+1}$, see Figure \ref{fig:schematic_process} for illustration. 

\subsubsection{Example: Burgers' Equation}
To illustrate the specifics of the time evolution, we consider Burgers' equation, with $f(u)=u^2/2$, approximated by $f_h=u_h^2/2$, where $u_h$ represents the degree $(2m+1)$ interpolant on either of the grids. We can write 
\begin{eqnarray}
(u_h)_t  +(f_h)_x &=& 0,  \nonumber\\
(u_h)_{tx}+(f_h)_{xx} &=& 0,\nonumber\\
(u_h)_{txx}+(f_h)_{xxx} &=& 0,\label{eq:discretized_PDE_1d}\\
  \vdots && \vdots \nonumber  
\end{eqnarray} 
where 
$$f_h = \sum_{k=0}^{2m+1}b_{k}(t)(x-x_{j-1/2})^{k}.$$
The coefficients $b_k$ are determined by truncated polynomial multiplication, that is $b_k = \frac{1}{2}\sum_{l=0}^kc_lc_{k-l}$. Insisting that the approximate solution $u_h$ satisfy equation (\ref{eq:discretized_PDE_1d}) at the cell centers $x=x_{j-1/2}$, we obtain
\begin{equation}
  \label{eq:discretized_system_1d}
  \begin{pmatrix}
    c'_{0}(t) \\
    c'_{1}(t) \\
    \vdots  \\
    c'_{2m}(t)\\
    c'_{2m+1}(t)
  \end{pmatrix}
  =
  \begin{pmatrix}
    b_1(t) \\
    2\;b_2(t) \\
    \vdots  \\
    (2m+1)\;b_{2m+1}(t)\\
    0
  \end{pmatrix}
  .
\end{equation}

While equation (\ref{eq:discretized_system_1d}) is valid for any cell, the initial data for each cell are different from one cell to another. For a more detailed explanation and open source implementations, see \cite{hagstrom2007experiments, icosahom2014,CHIDES}. 

  \setlength{\unitlength}{0.7mm}
  \begin{figure}[htb]
    \begin{center}
      \vspace{2.5cm}
      \begin{picture}(130,60)(-10,-30)
        \thicklines
        \matrixput(27.5,2.5)(50,0){2}(0,25){3}{\circle{2}}
        \matrixput(2.5,2.5)(50,0){3}(0,25){3}{\circle*{2}}
        \put(8,0){${\cal I}$}
        \put(15,0){$\rightarrow$}
        \put(58,0){${\cal I}$}
        \put(65,0){$\rightarrow$}
        \put(33,25){${\cal I}$}
        \put(40,25){$\rightarrow$}
        \put(83,25){${\cal I}$}
        \put(90,25){$\rightarrow$}
        \put(42,0){${\cal I}$}
        \put(35,0){$\leftarrow$}
        \put(92,0){${\cal I}$}
        \put(85,0){$\leftarrow$}
        \put(17,25){${\cal I}$}
        \put(10,25){$\leftarrow$}
        \put(67,25){${\cal I}$}
        \put(60,25){$\leftarrow$}
        \put(25,8){${\cal T}$}
        \put(26,15){$\uparrow$}
        \put(75,8){${\cal T}$}
        \put(76,15){$\uparrow$}
        \put(0,33){${\cal T}$}
        \put(1,40){$\uparrow$}
        \put(50,33){${\cal T}$}
        \put(51,40){$\uparrow$}
        \put(100,33){${\cal T}$}
        \put(101,40){$\uparrow$}
        \put(0,-15){$x_{j-1}$}
        \put(25,-15){$x_{j- \frac {1}{2}}$}
        \put(50,-15){$x_{j}$}
        \put(75,-15){$x_{j+ \frac {1}{2}}$}
        \put(100,-15){$x_{j+1}$}
        \put(115,0){$t_n$}
        \put(115,25){$t_{n + \frac {1}{2}}$}
        \put(115,50){$t_{n+1}$}
      \end{picture}
         \caption{Illustration of the numerical process in one dimensional Hermite Methods for a full time step. Solid circles represent the primal grid $\mathcal{G}^p$ and open circles represent the dual grid $\mathcal{G}^d$. ${\cal I}$ is the Hermite interpolation operator and ${\cal T}$ is the time evolution operator.}
         \label{fig:schematic_process}
    \end{center}
  \end{figure}

\subsection{Flux-Conservative Hermite Methods}
\label{sec:conservative_hermite}
The numerical flux $f_h$ obtained by the approach described above, is discontinuous at cell interfaces when the flux function $f(u)$ is nonlinear. Numerically, the discontinuity in the flux induces numerical stiffness. As a result, the time step often needs to be taken very small. To remedy this, we propose new flux-conservative Hermite methods that impose flux continuity by computing the numerical flux at cell edges, and then interpolate the numerical flux to cell center.

\begin{figure}[h]
  \begin{center} 
    \psfrag{xxx}[][][1][0]{$x$}
    \psfrag{yyy}[][][1][0]{$f(u)$}
   \includegraphics[width=0.6\textwidth]{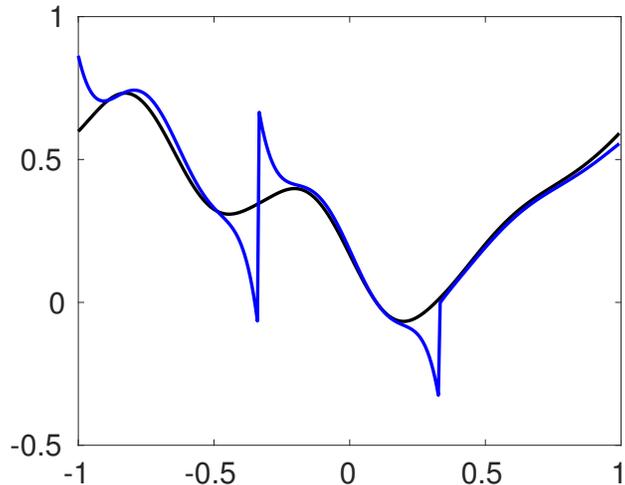}
  \caption{Numerical flux for Burgers' equation with random data, obtained by original Hermite (blue), flux-conservative Hermite (black). \label{fig:flux_random_data}}
  \end{center}
\end{figure}

To illustrate the difference between the original and flux-conservative Hermite schemes, we plot the numerical flux $f_h=u_h^2/2$ with $m=3$ and for $N_x=3$ cells in Figure \ref{fig:flux_random_data}. The numerical flux obtained using the original Hermite method, shown as the blue curve, has discontinuities at cell interfaces. On the other hand, the numerical flux obtained by the flux-conservative Hermite method, shown as the black curve, is continuous everywhere. 

\subsubsection{The Construction of the New Method}
The goal of our construction is to globally conserve the integral of $u_h$ and its $m$ first derivatives over one half time step, with $\Delta \hat{t} = \Delta t/2$.

In the flux-conservative methods, we assume that the solution on the primal grid at time $t_n$ is given by 
\begin{equation}
u^p(x,t_n) = \sum_{k=0}^{2m+1} c_{k}(t_n) \left(x-x_{j}\right)^k, x \in \mathcal{I}^d_j=[x_{j-1/2},x_{j+1/2}].
\end{equation}
Note that the degree of this polynomial is different than in the original method. We assume that the time stepping will be performed by an explicit one step method requiring stage values. The evolution will be carried out at the cell center where the stage values will be the derivative of the Hermite interpolant of the flux. As this interpolant is $m$ times differentiable at the edges, it will result in a conservative discretization. Now, the time evolution of the approximate solution entails the following steps.

{\it Computation of the stage fluxes at the cell edges:}  For simplicity, assume that we use the classic fourth order Runge-Kutta, then to construct the Hermite interpolants for the four stages we first compute 

\begin{align*}
    F_1^p & = f\left(u^p\right), \\
    F_2^p & = f\left(u^p+\frac{\Delta \hat{t}}{2}\frac{dF_1^p}{dx}\right), \\ 
    F_3^p & = f\left(u^p+\frac{\Delta \hat{t}}{2}\frac{dF_2^p}{dx}\right), \\
    F_4^p & = f\left(u^p+\Delta \hat{t}\frac{dF_3^p}{dx}\right). 
\end{align*}
Note that inside the argument of $f$, we keep all the coefficients of the polynomials up to degree $2m+1$, while the nonlinearity $f$ itself, which typically is a higher degree polynomial, is truncated to degree $2m+1$.

{\it Reconstruction of solution and fluxes by Hermite interpolation:} Next, we construct $u^d$ and $F_i^d$, $i=1,\ldots,4$, the global Hermite interpolation polynomials of degree $(2m+1)$ of the solution and the flux, respectively. Let $w^d$ represent $u^d$ or $F_i^d$ and $w^p$ represent $u^p$ or $F_i^p$. Then,
   \begin{equation}
     \label{eq:hermite_interpolant_conservative_1d}
     w^d(x,t_n) = \sum_{k=0}^{2m+1} d_{k}(t_n) \left(x-x_{j-1/2}\right)^k, \ \ x \in \mathcal{I}^p_j=[x_{j-1},x_{j}], 
   \end{equation} 
   where the coefficients $c_k(t_n)$ are uniquely determined by the interpolation conditions (at cell edges)
\begin{equation}
  \frac{\partial^k}{\partial x^k}w^d(x_i,t_n) = \frac{d^k}{dx^k}w^p(x_i,t_n), \ \ k=0,\ldots,m,\ \ i=j-1,j.
\end{equation}

{\it Time evolution:} Let the coefficients of $u^d$ be $c_k$ and the coefficients of $F_i^d$ be $b_k^{(i)}$. Then, again assuming RK4, we have that for $k=0,\ldots,2m+1$, 
\begin{equation*}
  c_k(t_n+\Delta \hat{t}) = c_k(t_n) + \frac{\Delta \hat{t}}{6} {(k+1)(b_{k+1}^{(1)}(t_n)+2b_{k+1}^{(2)}(t_n)+2b_{k+1}^{(3)}(t_n)+b_{k+1}^{(4)}(t_n))}. 
\end{equation*}
The updated solution on the dual grid is thus
\begin{equation*}
u^d(x,t_n+\Delta \hat{t}) = \sum_{k=0}^{2m+1} c_{k}(t_n+\Delta \hat{t}) \left(x-x_{j-1/2}\right)^k,\ \ x \in \mathcal{I}^p_j=[x_{j-1},x_{j}].
\end{equation*}
     
{\it Repeat on dual grid:} At the end of the half time step, we have evolved the degree $2m+1$ polynomial $u^d$. We then repeat the above process to obtain $u^p$ at time $t=t_{n+1}$.

We note that unlike the original method, the number of degrees of freedom that we keep is twice as large, representing an increase in memory requirement. The number of floating point operations, however, to leading order, is the same as for the original method (see the complexity analysis below). 

\subsubsection{Conservation}
We now consider the conservation properties of the above scheme. Due to the fact that the $s$ fluxes used in the stages have $m$ continuous derivatives we immediately find that for periodic problems the following conservation statements hold.  

{\bf Theorem:}
Assume we use the flux-conservative Hermite method to evolve $u_t+f(u)_x=0$ with periodic boundary conditions. Further assume that $u^d(t,x)$ is the periodic degree $2m+1$ Hermite interpolating polynomial and that $F^d_i$, $i=1,\ldots,s$ are the periodic degree $2m+1$ polynomials Hermite interpolating the fluxes. Further, let the coefficients $c_k(t)$ of $u^d$ on a cell be evolved from time $t = t_n$ to $t = t_n+\Delta \hat{t}$ by the one step method
  \[
  c_k(t_n+\Delta \hat{t}) - c_k(t_n) + \Delta \hat{t} \sum_{i=1}^s \alpha_i (k+1)b_{k+1}^{(i)}(t_n)=0, \ \  k = 0,\ldots,2m+1, 
  \]
  where $b_k^{(i)}$ are the coefficients of $F^d_i$. 
  Then, the following conservation statement holds.
  \begin{multline}
    \sum_j \int_{x_{j-1}}^{x_j} \frac{\partial^k}{\partial x^k} u^d(t_n+\Delta \hat{t},x) dx = \sum_j \int_{x_{j-1}}^{x_j} \frac{\partial^k}{\partial x^k} u^d(t_n,x) dx, \ \ k = 0,\ldots,2m+2-s. 
  \end{multline}

{\bf Proof:}
From the RK time stepping method for conservation laws (\ref{eq:conservation_PDE}) as
  \begin{equation}
    \label{eq:RK_cons_law_1D}
  \frac{u^d(t_{n+1/2})-u^d(t_n)}{\Delta \hat{t}}=\Delta \hat{t}\sum_i\alpha_i\frac{dF^d_i}{dx},
\end{equation}
together with the continuity of $m$ first derivatives of  each of the $F^d_i$, the result follows immediately from the update formula. Note also that $F^d_{i+1}$ is one order less accurate than $F^d_{i}$ due to flux differentiation during stage $i$. \\

\begin{figure}
    \begin{tikzpicture}
      \draw(-6,0) node{original};
      \draw[-](-4,0)  -- (0,0) -- (4,0);

      \draw(-4,0)node[circle,inner sep=2pt, fill=black,label={below:{$x_{j-1}$}}](xl){};
      \draw(0,0)node[circle,inner sep=2pt, fill=black,label={below:{$x_{j}$}}](xm){};
      \draw(4,0)node[circle,inner sep=2pt, fill=black,label={below:{$x_{j+1}$}}](xr){};
      \draw(-2,0)node[label={above:{$(\mathcal{I}u)_{j-\frac{1}{2}}$}}](il){};
      \draw(-2,0)node[rectangle, fill = black, label={below:{$f_{j-\frac{1}{2}}(\mathcal{I}u)$}}](fl){};
      \draw(2,0)node[label={above:{$(\mathcal{I}u)_{j+\frac{1}{2}}$}}](ir){};
      \draw(2,0)node[rectangle, fill = black, label={below:{$f_{j+\frac{1}{2}}(\mathcal{I}u)$}}](fr){};
      \draw[->](xl) to [bend left = 45](il);
      \draw[->](xm) to [bend left = 45](ir);
      \draw[->](xm) to [bend right = 45](il);
      \draw[->](xr) to [bend right = 45](ir);          

      \draw(-6,-2.5) node{flux conservative};
      \draw[-](-4,-2.5)  -- (0,-2.5) -- (4,-2.5);
      \draw(-4,-2.5)node[circle,inner sep=2pt, fill=black,label={below:{$x_{j-1}$}}](xl){};
      \draw(0,-2.5)node[circle,inner sep=2pt, fill=black,label={below:{$x_{j}$}}](xm){};
      \draw(4,-2.5)node[circle,inner sep=2pt, fill=black,label={below:{$x_{j+1}$}}](xr){};
      \draw(-4,-2.5)node[label={above:{$f_{j-1}(u)$}}](fl){};
      \draw(0,-2.5)node[label={above:{$f_j(u)$}}](fm){};
      \draw(4,-2.5)node[label={above:{$f_{j+1}(u)$}}](fr){};
      \draw(-2,-2.5)node[rectangle, fill = black, label={below:{$(\mathcal{I}f)_{j-\frac{1}{2}}$}}](fl){};
      \draw(2,-2.5)node[rectangle, fill = black, label={below:{$(\mathcal{I}f)_{j+\frac{1}{2}}$}}](fr){};
      \draw[->](xl) to [bend left = 45](fl);
      \draw[->](xm) to [bend left = 45](fr);
      \draw[->](xm) to [bend right = 45](fl);
      \draw[->](xr) to [bend right = 45](fr);
    \end{tikzpicture}
    \caption{Original vs. Flux Conservative Hermite Methods. Here, we dropped the subscript $h$ in all the computed quantities for compactness.}
    \label{fig:orig_vs_conservative_hermite}
\end{figure}
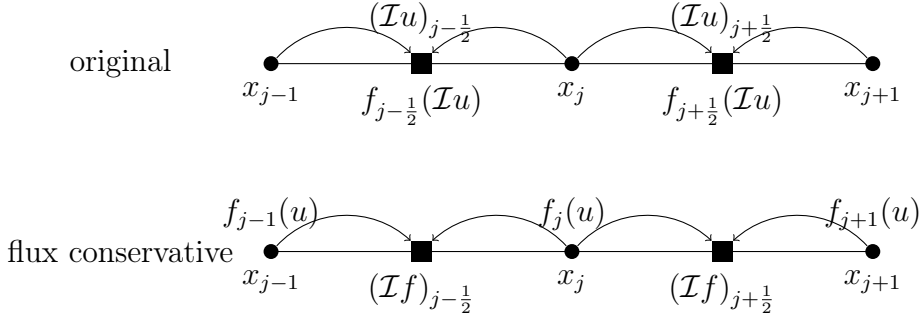

To summarize, in the original Hermite scheme, computation of numerical fluxes is performed at the cell center using the interpolated solution. The flux-conservative Hermite scheme requires a computation (and storage) of numerical fluxes at the cell edges and the interpolation of those fluxes to the cell center. Refer to Figure \ref{fig:orig_vs_conservative_hermite} for an illustrative comparison between the schemes.

\subsection{The Flux-Conservative Hermite Method in Two Dimensions} \label{sec:new_hermite_2D}
Now, let us consider a conservation law
\begin{equation}
  \label{eq:conservation_PDE_2D}
  u_t+(f(u))_x+(g(u))_y=0,
\end{equation}
on the domain $D = [x_L,x_R]\times[y_B,y_T]$. Let $\mathcal{G}^p$ and $\mathcal{G}^d$ be the primal and dual grids, defined as
\begin{equation}
\mathcal{G}^p = \{(x_{i},y_{j})\} = (x_L+ih_x,y_B+jh_y), \: i = 0,\ldots,N_x,\:j = 0,\ldots,N_y, 
\end{equation}

\begin{multline}
  \mathcal{G}^d =\left\{(x_{i+1/2},y_{j+1/2}) = \left(x_L+\left(i+\frac{1}{2}\right)h_x,y_B+\left(j+\frac{1}{2}\right)h_y\right)\right\}, \\
  i = 0,\ldots,N_x-1, \: j = 0,\ldots,N_y-1,
\end{multline}
where $h_x = (x_R-x_L)/N_x$ and $h_y = (y_T-y_B)/N_y$ are the distances between two adjacent nodes in $x$ and $y$ directions respectively.

The extension of the flux-conservative method from one dimension is straightforward. Writing equation (\ref{eq:conservation_PDE_2D}) as $u_t = -f(u)_x-g(u)_y$ and letting $u^d(t)$ represent the two-dimensional tensor product Hermite interpolant of the data on the primal grid we can write the RK4 time stepping of $u^d(t_n)$ to time $t = t_n + \Delta \hat{t}$ as 
\begin{equation}
  \label{eq:RK_cons_law_2D}
  \frac{u^d(t_n+\Delta \hat{t})-u^d(t_n)}{\Delta \hat{t}}=\frac{K_1^d+2K_2^d+2K_3^d+K_4^d}{6}.
\end{equation}
The left hand side of equation (\ref{eq:RK_cons_law_2D}) is an approximation to $u_t$ and the right hand side is an approximation to stage values of $-(f(u)_x+g(u)_y)$. Similar to the one dimensional case we have 
\begin{align*}
  K_1^d      & = -(F_1^d)_x-(G_1^d)_y, \\
  K_2^d      & = -(F_2^d)_x-(G_2^d)_y, \\
  K_3^d      & = -(F_3^d)_x-(G_3^d)_y, \\
  K_4^d      & = -(F_4^d)_x-(G_4^d)_y. 
\end{align*}
Here, for example, $F_i^d$ is the degree $2m+1$ tensor product polynomial that interpolates $f(u^p+\gamma_i \Delta \hat{t} F^p_i)$ and its $m$ first derivatives in $x$ and $y$ at the four adjacent primal grid-points.

\subsection{Comparison of Computational Costs}
The time evolution portion of the Hermite methods are performed by a one step method with $n_K$ stages, involving computation of the flux function, interpolation of the solution and, in the flux-conservative method, the fluxes, and differentiation of fluxes. For the purpose of this comparison, we assume Burgers' flux function $f(u) = u^2/2$ in 1D or $f(u) = g(u) = u^2/2$ in 2D. Each 1-dimensional Hermite interpolation is equivalent to a multiplication by a $(2m+2)$ by $(2m+2)$ matrix. If we use the recipe above, each 2-dimensional Hermite interpolation corresponds to ${2\times(2m+2)}$ one-dimensional interpolations. The factor $(2m+2)$ comes from the the fact that the $y$ dimension brings in $(2m+2)$ interpolations in 1D and the multiplicative factor $2$ comes from the fact that we interpolate in $y$ direction on both the left and right edges of the cell. In 3D, we have 4 interpolations in the $z$ direction, each brings in $(2m+2)$ times interpolations in 2D, and so on. We summarize the cost of the method, corresponding to the number of multiplications involved, in Table \ref{tab:cost_old_vs_new}. The number of interpolation in the flux-conservative Hermite method is $n_K d \;2^{2^{d-1}}$ more than the original Hermite method. We note that the flux-conservative scheme can be simplified to just two flux interpolations by adding up the $F$'s and $G$'s, but in this paper, we interpolate each flux separately. There is also an additional cost of differentiation at cell corners, but it is negligible compared to the cost of interpolation.
\begin{table}[]
  \begin{center}
    \begin{tabular}{|l|c|c|}
      \hline
                                    & {\bf flux computation}    & {\bf interpolation}      \\
      \hline
      original                      & $n_Kd(2m+2)^{2d}$       & $2^{2^{d-1}}(2m+2)^{d+1}$      \\ 
      \hline
      flux-conservative             & $n_Kd(2m+2)^{2d}$        & $(1+n_Kd)2^{2^{d-1}}(2m+2)^{d+1}$   \\ 
      \hline
    \end{tabular}
  \end{center}
  \caption{Comparison of Costs in original and flux-conservative Hermite methods, $n_K$ is the number of stages in Runge-Kutta method, $d$ is the spatial dimension. \label{tab:cost_old_vs_new}}
\end{table}

\section{The Entropy Viscosity Method}
\label{sec:EV_intro}
Given a PDE on the form (\ref{eq:conservation_PDE}), there exists an entropy function $E(u)$ and its corresponding entropy flux function $F(u)=\int E'(u)f'(u)\; du$ such that the entropy residual satisfies
\begin{equation}\label{eq:entropy_residual}
  r_{EV}(u) \equiv E_t(u) + \nabla \cdot F(u) \leq 0.
\end{equation}
This inequality can be used to select the physically correct solution to (\ref{eq:conservation_PDE}) or (\ref{eq:conservation_PDE_2D}). The direction of the inequality can vary from one problem to another, but the inequality becomes strict only at shocks. In essence, the entropy viscosity (EV) method uses the fact that the residual approaches a Dirac delta function centered at shocks to construct a nonlinear dissipation. The resulting dissipation is small away from shocks and just sufficient amount at a shock. The details of EV for conservation laws are described in detail in \cite{Guermond2011conservation} but we summarize its most important features here.

Consider the conservation law whose right hand side has been replaced by a viscous term, ${u_t+\nabla\cdot{f}(u) = \nabla \cdot (\nu \nabla u)}$, with $\nu=\nu_h(x,t)$ given by 
\begin{equation}\label{eq:final_viscosity}
\nu_h(x,t):= \min(\nu_{EV},\nu_{\max}),
\end{equation} 
where $\nu_{EV}$ is the entropy-based viscosity and $\nu_{\max}$ is a viscosity whose size depends on the largest eigenvalue of the flux function $f(u)$, representing the maximum wave speed. The discretized entropy-based viscosity $\nu_E$ is then given by 
\begin{equation}\label{eq:entropy_visc}
\nu_{EV}(x,t) = \alpha_{EV} C_1(u_h)h^{\beta} |r_{EV}(u_h)|,
\end{equation}
where $\beta$ is a positive scalar, $\alpha_{EV}$ is a user defined constant and $C_1(u_h)$ is some PDE-specific normalization.

At shocks, the entropy residual approaches Dirac delta function, so we instead use   
\begin{equation}\label{eq:visc_max}
  \nu_{\max}(x,t)=\alpha_{\max}\; h \max\limits_{y \in V_x} C_2(u_h(y,t)) ,
\end{equation}
where $\alpha_{\max}$ is another user defined constant, $C_2$ serves as the maximum wave speed and $V_x$ is some neighborhood of $x$. The $V_x$ neighborhood can either be ``local'', i.e. containing only a few cells around $x$, or ``global'', i.e. $V_x=D$, where $D$ is the domain of the PDE. In this work, we use global $V_x$.

In recent papers, the parameter $\beta$ is chosen to be 2, but we found that this may prevent convergence for moving shocks, see \cite{icosahom2016} where we also argue that $\beta=1$ is a more appropriate choice. In essence our argument is as follows. As the entropy residual approaches a Dirac delta distribution, a consistent discretization of the residual with a single shock must satisfy  
\begin{equation}
\sum_{j=0}^{N_x-1}h_j (r_{EV})_j = 1,
\end{equation}
where $h_j = x_{j+1}-x_j$. Thus, we expect $(r_{EV})_j\sim h_j^{-1}$ near the shock. When $\beta=2$, the two terms $\nu_{EV}$ and $\nu_{\max}$ are both $\mathcal{O}(h)$. Since the parameters are tuned on a coarse grid and the terms $C_1$ and $C_2$ in (\ref{eq:entropy_visc})-(\ref{eq:visc_max}) also changes with the grid size, the selection of the minimum of $\nu_{EV}$ and $\nu_{\max}$ does not necessarily  ``converge'' as the grid gets refined. If instead we choose $\beta=1$, then $\nu_{EV} = \mathcal{O}(1)$ while $\nu_{\rm max} = \mathcal{O}(h_j)$, and the particular choice of $ {\alpha}_{EV}$ is thus irrelevant in the limit $ h_j \to 0 $ as the selection mechanism will eventually select $\nu_{\rm max}$ at the shocks.

While the explicit formula for $C_1$ and $C_2$ varies from one PDE to another, the core of the entropy viscosity method remains the same. The size of the entropy residual gives us a sense of relative distance to the shock, which is then used to take the following actions:   near a shock, EV uses sufficiently large dissipation, $\nu_h = \nu_{\max}$, and away from a shock, EV uses entropy-based dissipation, $\nu_h = \nu_{EV}$.

\section{Numerical Examples}\label{sec:numerical_exp}
We start by confirming that the rates of convergence, $2m+1$, (in space) are the same for the new flux-conservative method as for the original method.  

\subsection{Convergence for a Smooth Solution}
We solve Burgers' equation on the domain $x\in [-\pi, \pi]$ and impose periodic boundary conditions. The initial data  is $u(x,0) = - \sin (x) + 0.3$ and we evolve the solution until time $t=0.4$. The timestep is chosen as $\Delta t = {\rm CFL} \, \,  h_x /$ $\max_x |u(x,0)|$, with ${\rm CFL} = 0.1$.
\begin{table}[]
\caption{Convergence study of smooth solution to Burgers' equation \label{tab:burger_smooth}}
\begin{center}
\begin{tabular}{| l | r | r | r | r | r | r | r |}
\hline
$h_x$ & $\pi / 2 $ & $\pi / 4 $ & $\pi / 8 $ & $\pi / 16 $ & $\pi / 32$  \\
\hline
\hline
$l_\infty$-err $m=1 $ &2.30(-01) &  5.85(-02) & 1.09(-02) &  1.42(-03)  & 1.80(-04) \\ 
Rate & & 1.97  &   2.43 &  2.94 & 2.98 \\ 
\hline
$l_\infty$-err $m=2 $ & 4.85(-02) &  5.47(-03) &  2.19(-04) &  7.25(-06)  & 1.97(-07)  \\
Rate && 3.15 &   4.64 &    4.92  &  5.20  \\ 
\hline
$l_\infty$-err $m=3 $ &1.09(-02)  & 6.59(-04) &  7.71(-06) &  4.70(-08) &  2.73(-10) \\
Rate && 4.04 & 6.42 & 7.36 & 7.43 \\
\hline
\end{tabular}
\end{center}
\end{table}%

In Table \ref{tab:burger_smooth}, we display the maximum error at the final time computed against a reference solution computed using $m = 5$ and $h_x = \pi / 64$. As can be seen the rates of convergence appear to approach the predicted rate $2m+1$.

We next present a sequence of experiments displaying the performance of the Hermite-Runge-Kutta-4-Entropy-Viscosity method for Euler's equations (with artificial viscosity). 

\subsection{Euler's equations in One Dimension}
We consider Euler's equations which represent conservation of mass, momentum, and energy, 

\begin{equation}
  \label{eq:Euler_PDE}
  \begin{pmatrix}
    \rho\\
    \rho u\\
    E
  \end{pmatrix}
  _t
  +
  \begin{pmatrix}
    \rho u\\
    \rho u^2+p\\
    (E+p)u
  \end{pmatrix}
  _x = 
  \begin{pmatrix}
    0\\
    0\\
    0
  \end{pmatrix}.
\end{equation}
Here, $\rho$ is the density, $\rho u$ is the momentum, $u$ is the velocity, and $E$ is the energy.
Furthermore, we assume an ideal gas with the equation of state
\begin{equation}
  \label{eq:ideal_gas_law_1D}
  E = \frac{p}{\gamma-1}+\frac{\rho u^2}{2},
\end{equation} 
where $\gamma=1.4$ is the adiabatic index and $p$ is the pressure. 

To regularize equation (\ref{eq:Euler_PDE}), we add a viscous term $(\nu(\rho,\rho u, E)_x)_{x}$, where the coefficient $\nu$ is obtained using the entropy viscosity method. Thus, the viscous Eulers' equations can be written as
\begin{equation}
  \label{eq:Euler_PDE_visc}
  \begin{pmatrix}
    \rho\\
    \rho u\\
    E
  \end{pmatrix}
  _t
  +
  \begin{pmatrix}
    \rho u - \nu \rho_x\\
    \rho u^2+p - \nu (\rho u)_x\\
    (E+p)u - \nu E_x
  \end{pmatrix}
  _x = 
  \begin{pmatrix}
    0\\
    0\\
    0
  \end{pmatrix}.
\end{equation}
We note that an alternative to this simple viscosity would be to use the full Navier-Stokes equations.  

\subsubsection{Entropy Viscosity (EV) method for 1D Euler's equations}
\label{sec:EV_1D}
The discretized viscosity coefficient $\nu=\nu_h$ is given in terms of primitive variables $\rho,p,\text{and }u$,  
\begin{align}
  \label{eq:EV_Euler_1D}
     \nu_h     =& \min(\nu_{\max},\nu_{EV}),\\
     \nu_{EV}   =& \alpha_{EV} \: h_x \: \rho_h(x,t) |r_{EV}(x,t)|,\\ 
     \nu_{\max} =& \alpha_{\max} \:h_x \: \rho_h(x,t) \max \limits_{y\in D} \left(|u_h(y,t)|+ \sqrt{\gamma T_h(y,t)}\right),
  \end{align}
  where $T_h=p_h/\rho_h$ is the temperature, $h_x$ is the grid size and
  \begin{equation}\label{eq:Euler_residual_1D}
    r_{EV} = \partial_tS_h+((uS)_h)_x \geq 0,
  \end{equation}
  is the entropy residual for the entropy function $S_h(p_h,\rho_h) = \frac{\rho_h}{\gamma-1}\log\left(\frac{p_h}{\rho_h^{\gamma}}\right)$ and its corresponding entropy flux $(uS)_h$. 

\subsubsection{An Improved Entropy Viscosity}
\begin{figure}[ht!]
    \begin{center}
      \includegraphics[width=0.45\textwidth]{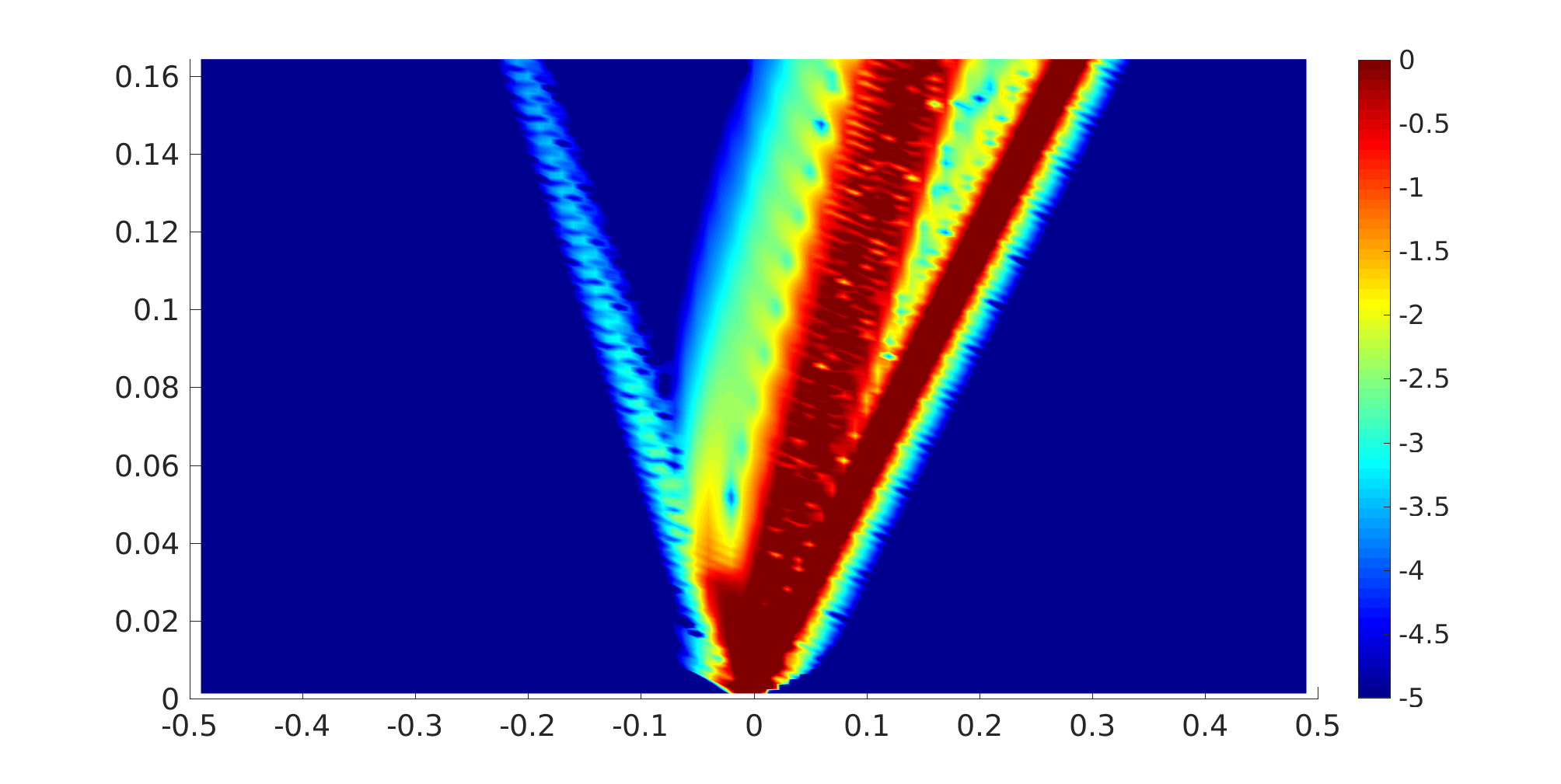} \hspace{0.2cm}
      \includegraphics[width=0.45\textwidth]{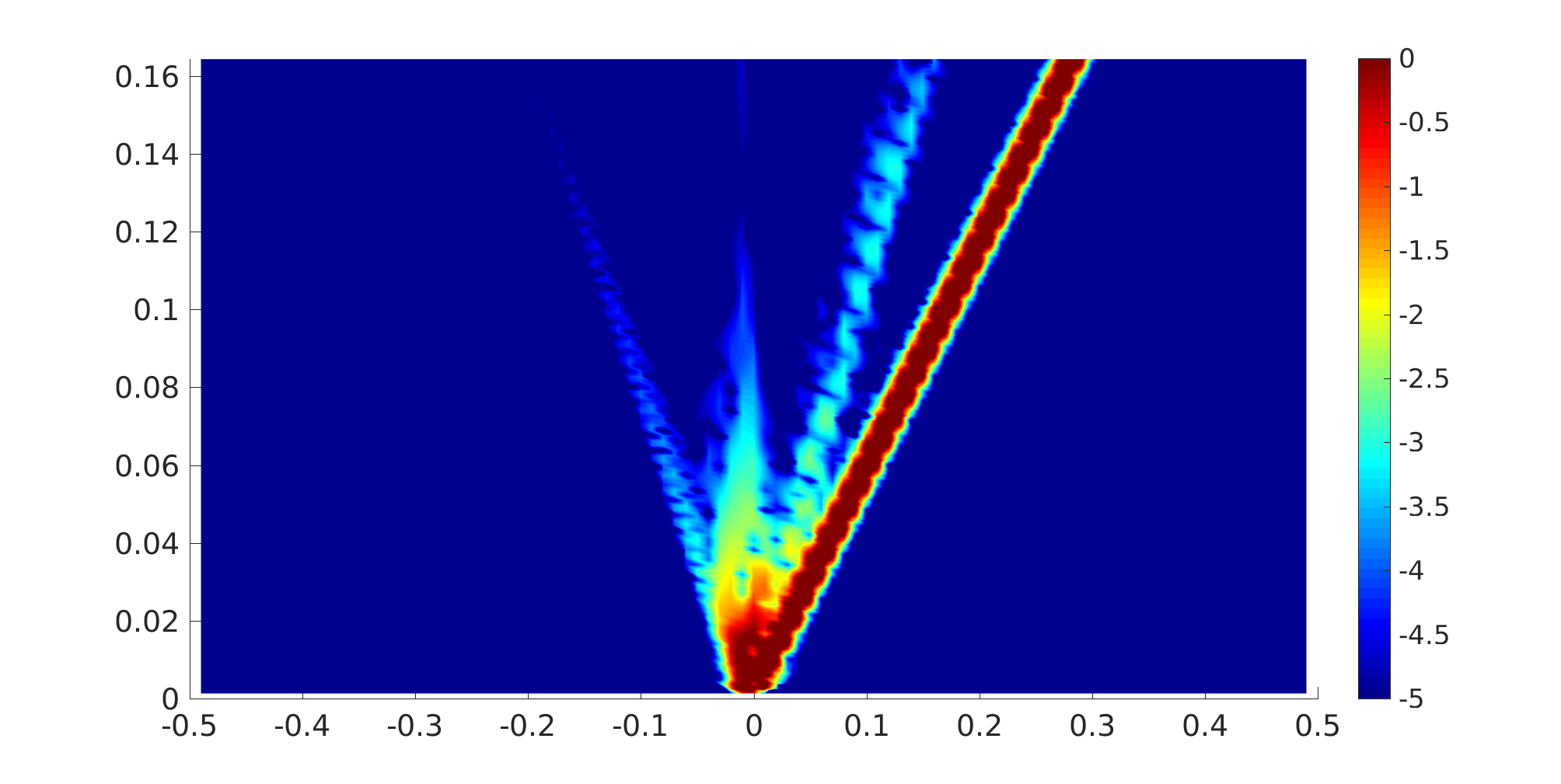}
    \caption{Space time diagram of the magnitude of entropy residual $|r_E|$ (left) and $|\Delta u\,r_E|$ (right) for Sod's problem. Blue is small, red is big. Simulations performed with $\nu_S \propto r_E$ (left) and $\nu_S \propto \Delta u\, r_E$ (right). With the new sensor $|\Delta u\,r_E|$, the residual, hence the viscosity is driven to zero along contact discontinuity (thicker red line in the middle disappears with the new sensor). \label{fig:res_xt_diag}}
\end{center}
\end{figure}

The entropy viscosity method discretizes the entropy residual using the numerical solution. In theory, the entropy residual is large at shocks, and zero at contact discontinuities (where no artificial viscosity is needed). However, our experience is that the discretization of the entropy equation may also trigger the maximum viscosity at contact discontinuities. To the left in Figure \ref{fig:res_xt_diag}, we see a space-time diagram of the entropy residual for Sod's problem in logarithm scale. Note that a relatively large amount of residual is produced at the contact discontinuity.
\begin{figure}[ht!]
    \begin{center}
      \includegraphics[width=0.45\textwidth]{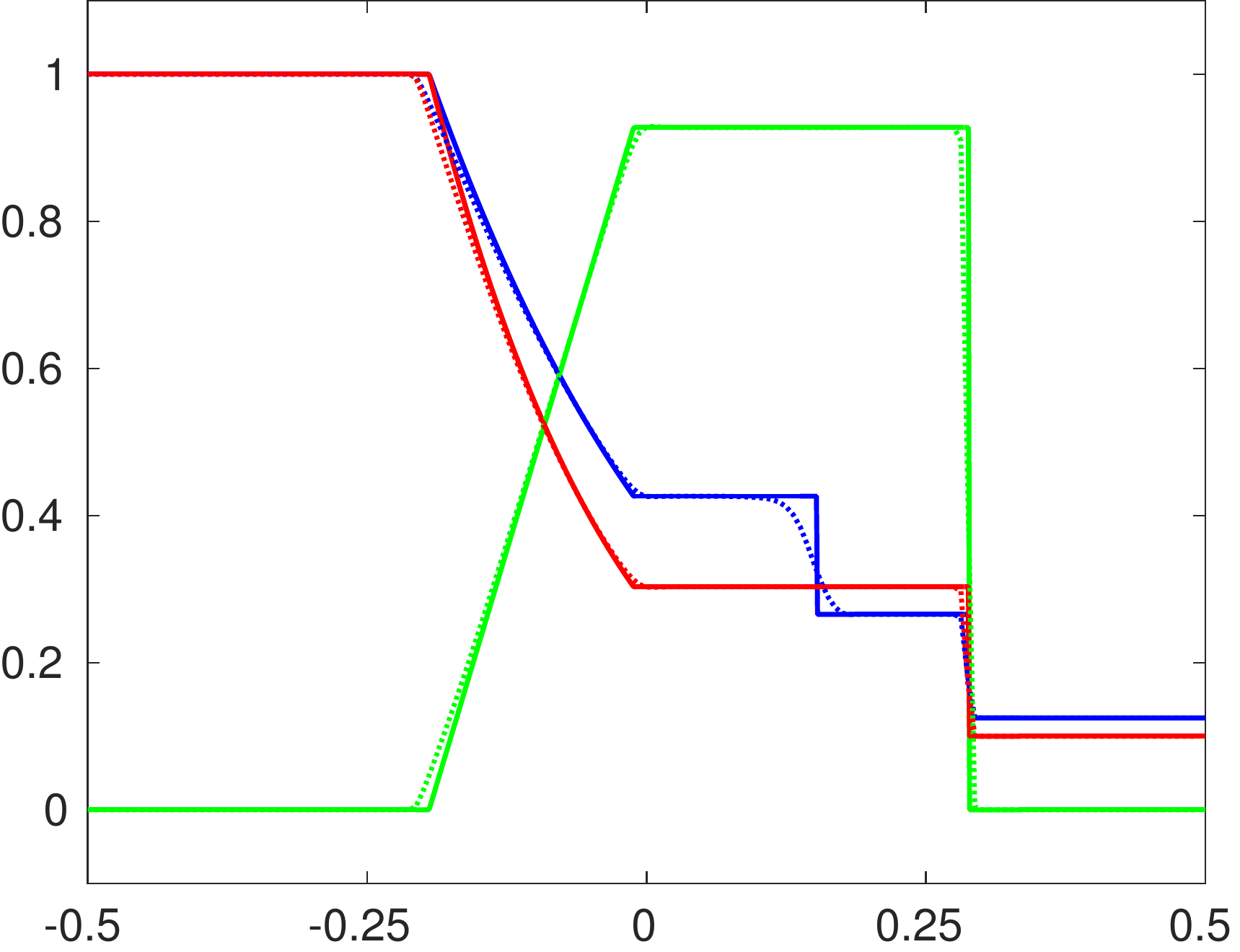} \hspace{0.2cm}
      \includegraphics[width=0.45\textwidth]{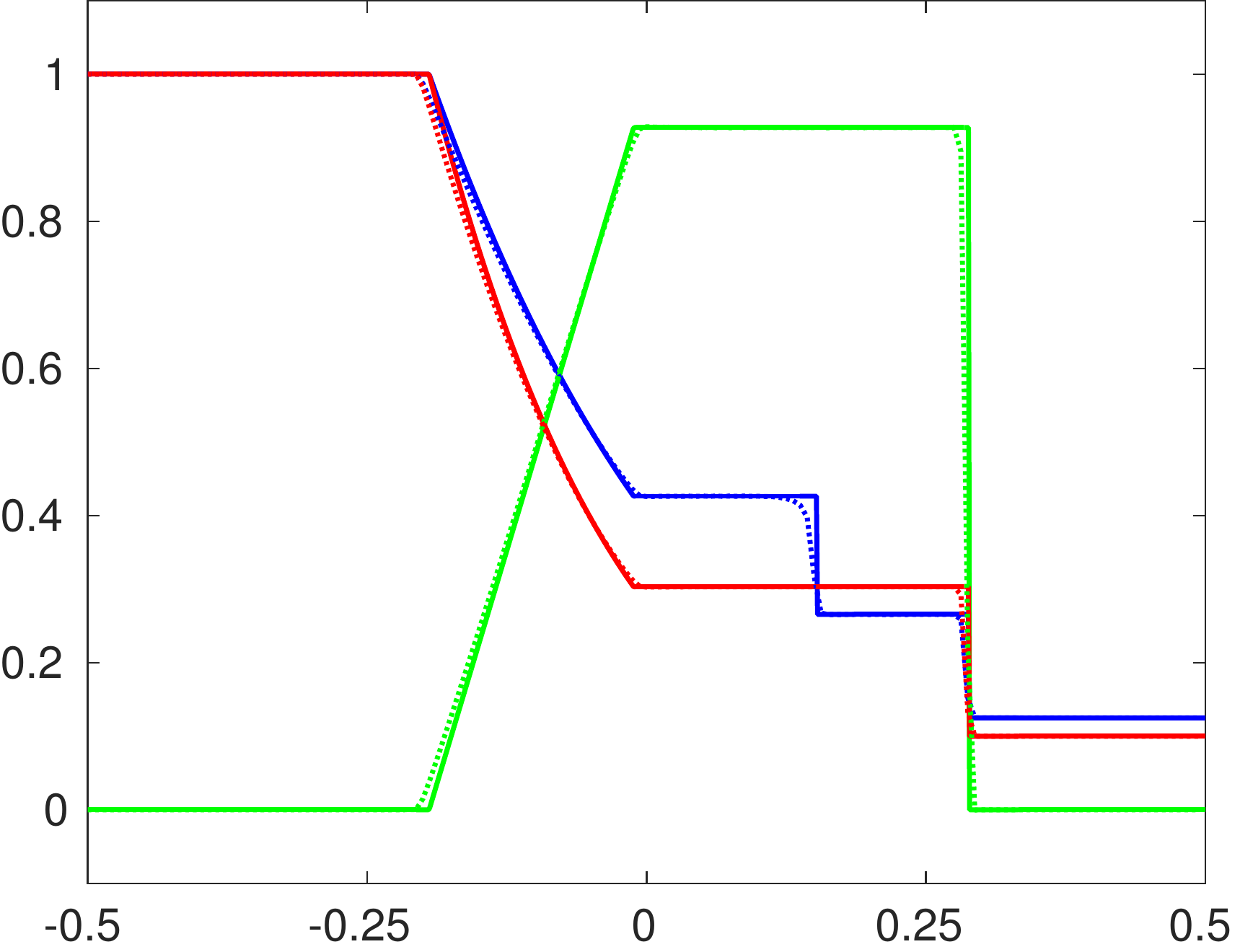}
    \caption{Numerical solution (dotted lines) obtained using entropy viscosity proportional to $|r_E|$ (left) and $|\Delta u\,r_E|$ (right) for Sod's problem are plotted against the exact solution (solid lines). We plot density (blue), velocity (green), and pressure (red). Contact discontinuity is sharper using the new sensor $|\Delta u\,r_E|$. \label{fig:sod_comparison}}
\end{center}
\end{figure}

To eliminate this undesired behavior along the contact discontinuity, we use the fact that  the velocity of the fluid, $u$, is a Riemann invariant along the second characteristic field. Since $u$ is smooth at the contact discontinuity but not at a shock, the product of $\Delta u_j$ and $(r_E)_j$ is small at contact discontinuities but still large at shocks. We incorporate the term $\Delta u$ into the improved entropy viscosity 
$$\nu_{EV}   = \alpha_{EV} \: h_x \: \rho_h(x,t) |\Delta u| |r_{EV}(x,t)|.$$

To this end we take $\nu_E$ and $\nu_{max}$ as a piecewise constant function on each cell. Thus, we compute the discretized density $\rho_h$, velocity $u_h$, temperature $T_h$, entropy function $S_h$ and entropy flux $uS_h$ at cell center in pointwise manner. Now, to get the entropy residual given in (\ref{eq:Euler_residual_1D}), we compute temporal and spatial derivatives using finite differences. Using the notation $S_h = S_j^n$ to denote the approximate flux function $S$ at $x=x_j$, $t=t_n$, we discretize  the term $\partial_tS_j^n$ using second order Backward Difference formula
\begin{equation}
  \partial_tS_j^n = \frac{3S_j^n-4S_j^{n-1}+S_j^{n-2}}{2\Delta t}.
\end{equation}
Similarly, the term $\partial_x (uS)_j^n$ is approximated by the centered finite difference 
\begin{equation}
  \label{eq:entropy_flux_discretization_1D}
\partial_x(uS)_{j}^n = \frac{(uS)_{j+1}^n-(uS)_{j-1}^{n}}{2h_x}.
\end{equation}

\subsection{Experiments in One Dimension with Euler's Equations}
We now present results obtained using our Hermite-RK4-EV method for a stationary shock, the Lax, the Sod, and the Shu-Osher problem. For experiments where we use more than one resolution, the EV parameters are tuned on the coarsest grid. In the 1D Euler's equations simulations, the timestep is chosen as $\Delta t = {\rm CFL} \,  h_x / \max_x |(u\pm c)(x,0)|$, where $c=\sqrt{\gamma p/ \rho}$ is the speed of sound, with ${\rm CFL}$ values given in Table \ref{tab:euler_conv}. 
\begin{table}[]
  \begin{center}
    \begin{tabular}{|l|c|c|c|c|c|c|}
      \hline
      $h_x$          & 1.25(-02)    & 6.25(-03)   & 3.13(-03)  & 1.56(-03) &  7.81(-04)  & 3.91(-04) \\
      \hline
      $L_1$ error    & 3.44(-03)   & 2.29(-03)   & 1.40(-03)  & 7.96(-04) &  4.46(-04)  & 2.49(-04)\\ 
      \hline
      Rate           &            & 0.59        & 0.71      & 0.81     &  0.83      & 0.84\\
      \hline
    \end{tabular}
  \end{center}
  \caption{Convergence study on Euler's equation with stationary shock. \label{tab:euler_conv}}
\end{table} 

\subsubsection{Stationary shock problem}
By solving the Riemann problem, we decide the states corresponding to a stationary shock. The goal of this experiment is to investigate the stability and accuracy of EV in the presence of shocks. Since small oscillations coming from shocks could potentially pollute the ``smooth'' regions, this is also a test for how well EV removes numerical artifacts.  The computational domain is $D=[-0.5,0.5]$ with the stationary shock given by  
\begin{equation}
    (\rho,u, p)(x,t) =
  \begin{cases}
    (0.84,1.08,0.56)      & x < 0, \\
    (1,0.9,0.71)          & x > 0. \\
  \end{cases} 
\end{equation}  
The boundary condition are imposed by setting the solution at the boundary so that it coincides with the solution at initial time. We perform a grid refinement study and report the errors in the density in Table \ref{tab:euler_conv}. We also present the ratio between successive errors. It appears that the rate of convergence for the $L_1$ error is approaching 7/8.

\subsubsection{Lax's and Sod's Shock Tube Problems}
Lax's and Sod's problems come from physical experiments in which a gas tube is separated by a membrane into two sections. The gas in each section is uniform in the $y$ and $z$ direction, so the problem is modeled as a 1-dimensional shock tube. The gas in the left section is kept at a different state than the gas in the right section. At time $t=0$, the membrane is punctured. In the problem setup, the Euler's equations are solved on the domain $D=[-0.5,0.5]$ with initial data
\begin{equation}
    (\rho,u, p)(x,0) =
  \begin{cases}
    (0.445,0.698,3.528) & x<0 \\
    (0.5,0,0.571)       & x>0 \\
  \end{cases} 
\end{equation}  
for Lax, and 
\begin{equation}
    (\rho,u,p)(x,0) =
  \begin{cases}
    (1,0,1)           & x<0 \\
    (0.125,0,0.1)     & x>0 \\
  \end{cases} 
\end{equation}  
for Sod. For both problems, we impose fixed boundary condition so that the solution on the boundary is the same as at the initial time. The solution is computed up to time $t=0.16$ for Lax's problem and time $t=0.1644$ for Sod's problem. 

The solution to Riemann problems such as Lax's and Sod's shock tubes contains 3 waves propagating from the discontinuity at the initial time. The second wave is a contact discontinuity, where the discontinuity is translated over time. The first and third waves are nonlinear, and can take either rarefaction waves or shock waves. 

\begin{figure}[htb]
    \includegraphics[width = 0.45\textwidth]{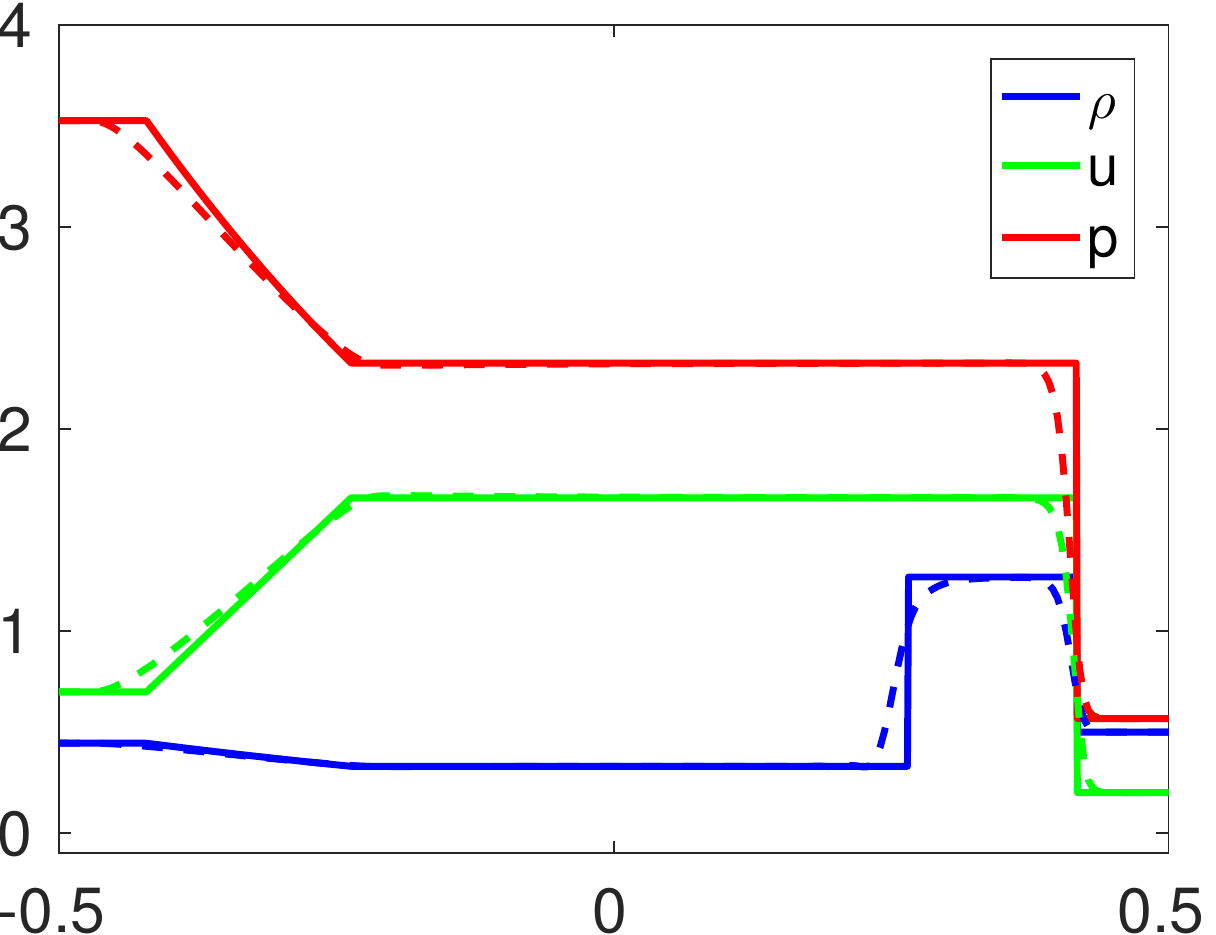}\hspace{0.5cm}
    \includegraphics[width = 0.45\textwidth]{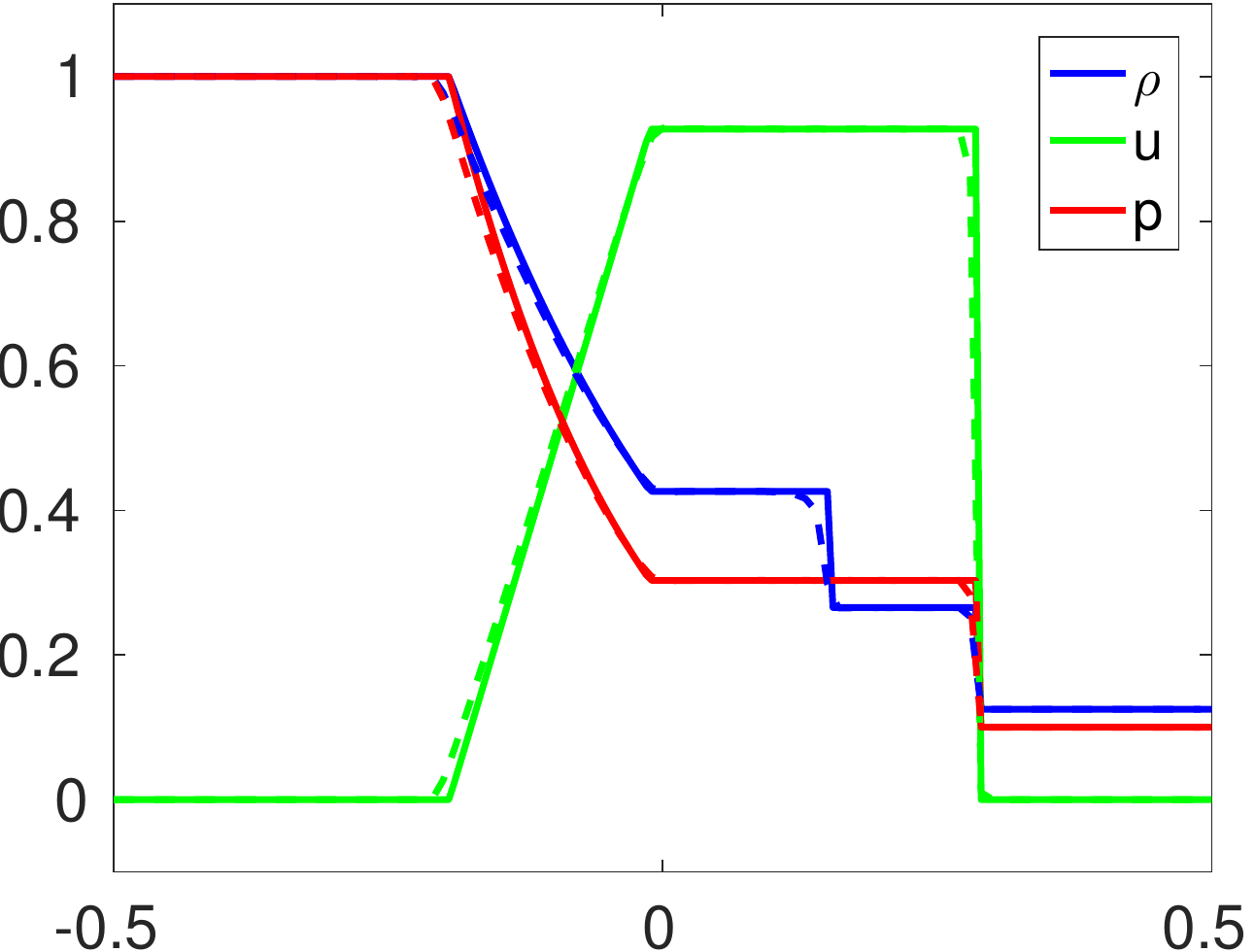}
    \caption{To the left: Lax shock tube, to the right: Sod shock tube. Dashed lines are the numerical solutions, solid lines are the exact solutions. Numerical solutions are obtained using $N_x=100$ elements. The color blue represents density, green represents velocity, and red represents pressure.}
    \label{fig:lax_sod}
  \end{figure}

The results for density $\rho$, velocity $u$ and pressure $p$ are plotted against the exact solution in Figure \ref{fig:lax_sod}. In each plot, we use $N_x = 100$ elements. The entropy viscosity parameters used are given in Table \ref{tab:parameters_1D}. In both problems, the shocks are resolved better than the contact discontinuities. Although the shock strength is only of medium size for both problems, some experts considered these tough test cases for non-characteristic-based high order schemes \cite{Shu1988439}. 

\begin{table}[]
  \begin{center}
    \begin{tabular}{|l|c|c|c|c|}
      \hline
      Problem      & {\rm CFL} & $m$  & $\alpha_{EV}$ & $\alpha_{max}$  \\
      \hline
      Lax          & 0.2 & 3  & 0.5 & 0.08  \\ 
      \hline
      Sod          & 0.15 & 3   & 0.2  & 0.08  \\ 
      \hline
      Shu-Osher    & 0.15 & 3  & 0.01  & 0.05  \\ 
      \hline
      Stationary shock & 0.2 & 3 & 10 & 0.3 \\ 
      \hline
    \end{tabular}
  \end{center}
  \caption{Parameters for examples in 1D Euler's  equations. \label{tab:parameters_1D}}
\end{table} 

\subsubsection{The Shu-Osher Problem}
\begin{figure}[htb]
    \centering
    \includegraphics[width = 0.5\textwidth]{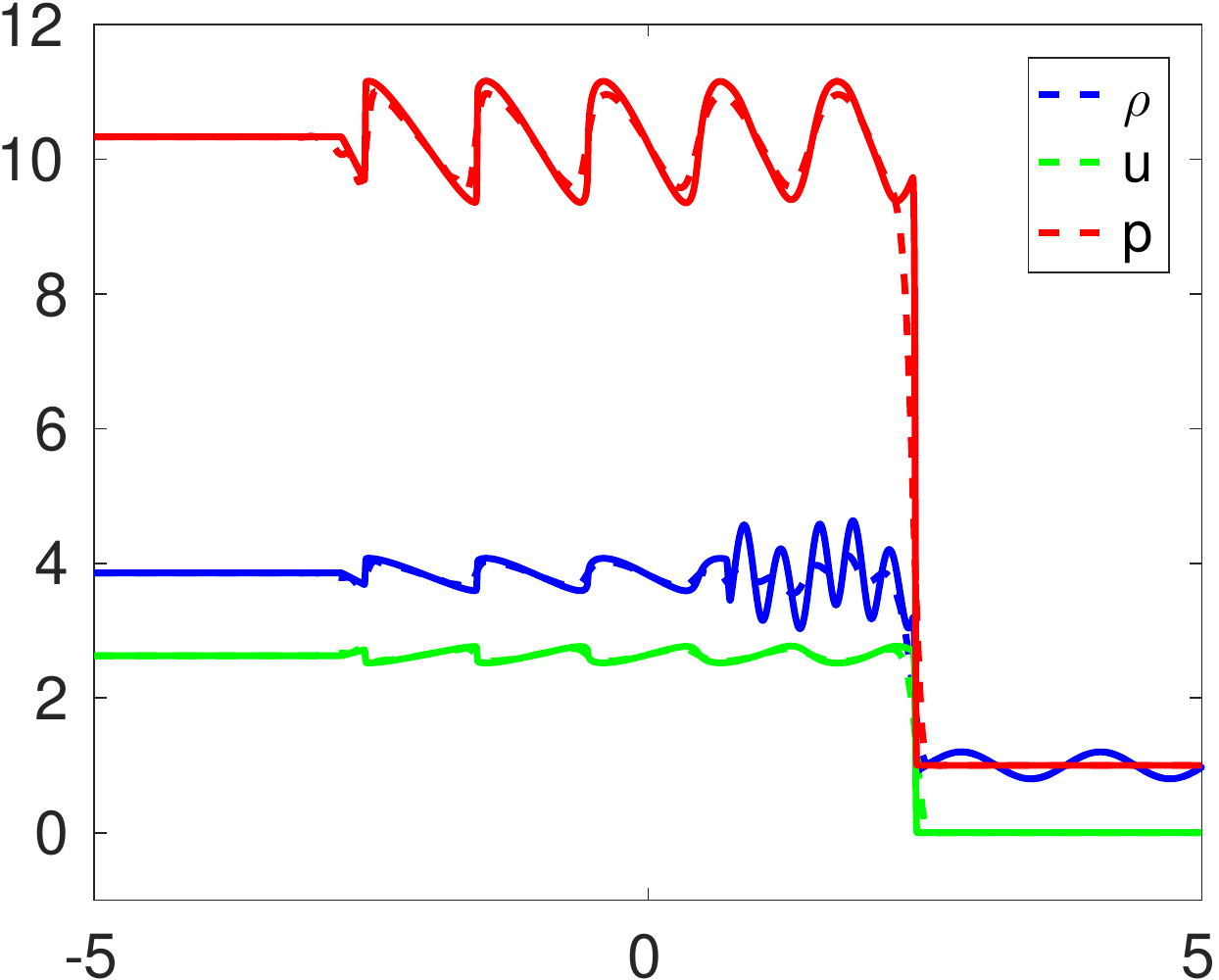}
    \caption{Shu-Osher problem. Dashed lines are the numerical solutions, solid lines are the ``exact'' solutions. Numerical solutions are computed using $N_x=80$ elements, ``exact'' solutions are computed on a much finer grid, with $N_x=1280$ elements. The color blue represents density, green represents velocity, and red represents pressure.}
    \label{fig:shu_osher}
  \end{figure}

The Shu-Osher problem poses difficulties for numerical methods due to the sinusoidal interacting with the shock. Here the domain is $D=[-5,5]$ with initial data
\begin{equation}
    (\rho,u, p)(x,0) =
  \begin{cases}
    (3.86,2.63,10.33)      & x<-4, \\
    (1+0.2\sin (5x),0,1)   & x>-4, \\
  \end{cases} 
\end{equation}
and with fixed boundary condition so that the solution on the boundary coincides with the solution at the initial time. The solution is computed up to time $t=1.8$ and compared against a computed solution on a much finer grid. We use $N_x=80$ to obtain the numerical solution in Figure \ref{fig:shu_osher}, where we interpolate the solution on to a finer grid. The ``exact'' solution is computed on a grid with $N_x=1280$. Note that even if we use a coarse grid, we can still get roughly the shape of the solution, especially away from the shock. However, when smooth waves are present (see blue oscillatory line to the left of shock) and too close to the shock, these waves get damped. 
  
\subsection{Euler's equations in Two Dimensions}
The two dimensional viscous Euler equations are given by 
\begin{equation}
  \label{eq:Euler_PDE_2D}
  \begin{pmatrix}
    \rho\\
    \rho u\\
    \rho v\\
    E
  \end{pmatrix}
  _t
  +
  \begin{pmatrix}
    \rho u\\
    \rho u^2+p\\
    \rho uv\\
    (E+p)u
  \end{pmatrix}
  _x
  +
  \begin{pmatrix}
    \rho v\\
    \rho uv\\
    \rho v^2+p\\
    (E+p)v
  \end{pmatrix}
  _y 
  = 
  \begin{pmatrix}
    0\\
    0\\
    0\\
    0
  \end{pmatrix}.
\end{equation}
Here, $\rho$ is the density, $\rho u$ and $\rho v$ are the momentum, $u$ and $v$ are the velocity in $x$ and $y$ direction respectively and $E$ is the energy.
Furthermore, we assume an ideal gas with equation of state
\begin{equation}
  \label{eq:ideal_gas_law_2D}
  E = \frac{p}{\gamma-1}+\frac{\rho (u^2+v^2)}{2},
\end{equation} 
where $\gamma=1.4$ is the adiabatic index and $p$ is the pressure. For all experiments below, the timestep is chosen as ${\Delta t = {\rm CFL} \, \,  h_x / \max_x |(u\pm c)(x,0)|}$, where $c=\sqrt{\gamma p/ \rho}$ is the speed of sound, with ${\rm CFL}$ values given in Table \ref{tab:parameters_2D}.

\begin{table}[]
  \begin{center}
    \begin{tabular}{|l|c|c|c|c|}
      \hline
      Problem      & $CFL$ & $m$   & $\alpha_{EV}$ & $\alpha_{max}$  \\
      \hline
      Explosion/implosion          & 0.2   & 3  & 0.1 & 0.2  \\ 
      \hline
      Vortex-shock interaction 1   & 0.2  & 3   & 0.01  & 0.04  \\
      \hline
      Vortex-shock interaction 2   & 0.2  & 3   & 0.05  & 0.07  \\ 
      \hline
      Jet                          & 0.2  & 3  & 0.03  & 0.2  \\ 
      \hline
    \end{tabular}
  \end{center}
  \caption{Parameters for examples in 2D Euler's  equations. \label{tab:parameters_2D}}
\end{table}

\subsection{Entropy Viscosity Method for Euler's Equations in Two Dimensions}
\label{sec:EV_2D}
The entropy viscosity is identical to the 1D version given in (\ref{eq:EV_Euler_1D}), with the exception that it takes the velocity in both directions into account.
\begin{align}
  \label{eq:EV_Euler_2D}
     \nu_h     =& \min(\nu_{\max},\nu_{EV}),\\
     \nu_{EV}   =& \alpha_{EV} \: h \: \rho_h(x,t) |r_{EV}(x,t)|,\\ 
     \nu_{\max} =& \alpha_{\max} \:h \: \rho_h(x,t) \max \limits_{y\in D} \left(\sqrt{u^2_h(y,t)+v^2_h(y,t)}+ \sqrt{\gamma T_h(y,t)}\right),
  \end{align}
  where $T_h=p_h/\rho_h$ is the temperature, $h=\min(h_x,h_y)$ is the grid size and
  \begin{equation}\label{eq:Euler_residual_2D}
    r_{EV} = \partial_tS_h+((uS)_h)_x + ((vS)_h)_y\geq 0.
  \end{equation}

  To discretize the entropy residual $r_{EV}$, we again use BDF for the time derivative and centered finite differences for the spatial derivatives,
  \begin{align}
    \label{eq:entropy_flux_discretization_2D}
    \partial_x(uS)_{jk}^n & = \frac{(uS)_{j+1,k}^n-(uS)_{j-1,k}^{n}}{2h_x}, \\
    \partial_y(vS)_{jk}^n & = \frac{(vS)_{j,k+1}^n-(vS)_{j,k-1}^{n}}{2h_y}. 
  \end{align}
  On the domain $[x_L,x_R]\times [y_B,y_T]$, we use the subscript $jk$ to indicate that the variable attached is evaluated at $x=x_L+jh_x$ and $y=y_B+kh_y$. 
  
\subsubsection{Explosion/Implosion Test}
First we solve a radially symmetric Riemann problem from Toro \cite{toro2013riemann}. The  computational domain is $D=[-1,1]\times[-1,1]$, and the initial data corresponding to an expanding wave is
\begin{equation}
    (\rho,u,v,p)(r,t) =
  \begin{cases}
    (1,0,0,1)      & r<0.4, \\
    (1,0,0,0.1)    & r>0.4. \\
  \end{cases} 
\end{equation}
For an imploding wave, the initial data is,
\begin{equation}
    (\rho,u,v,p)(r,t) =
  \begin{cases}
    (1,0,0,1)      & r>0.4, \\
    (1,0,0,0.1)    & r<0.4. \\
  \end{cases} 
\end{equation}

\begin{figure}
    \begin{center}
      \includegraphics[width = 0.45\textwidth]{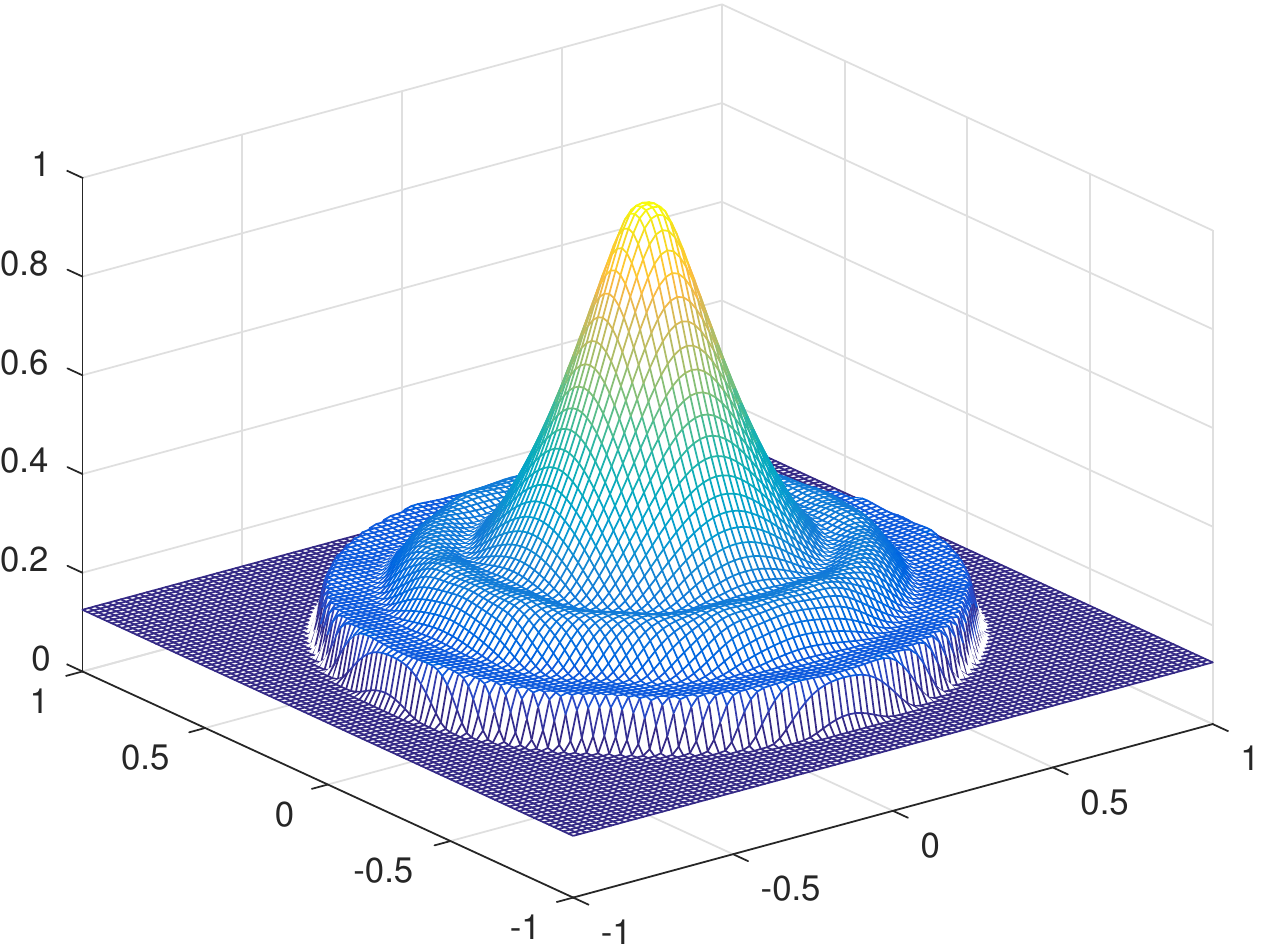}\hspace{0.5cm}
      \includegraphics[width = 0.45\textwidth]{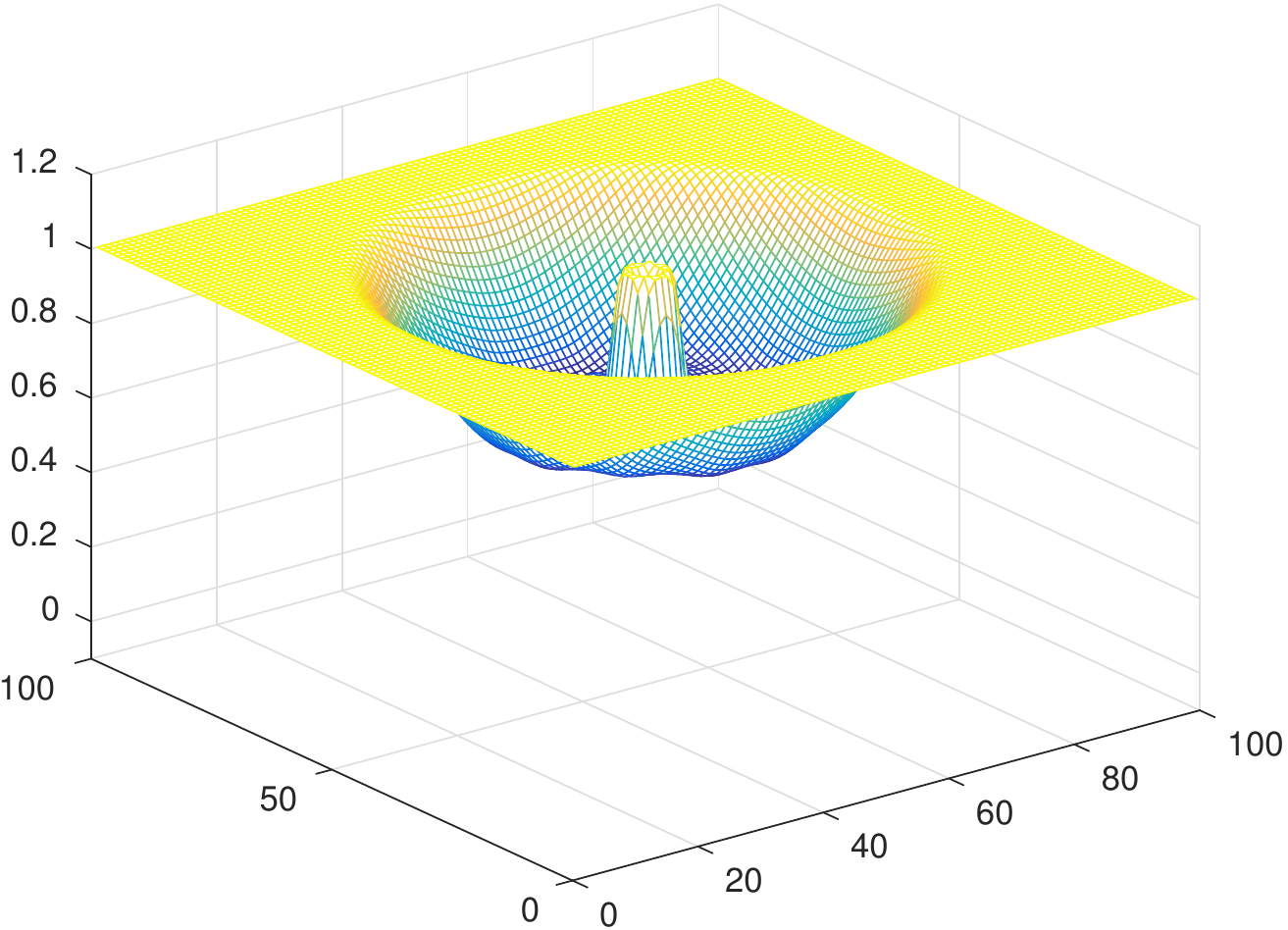}
    \caption{Solution to explosion problem at time $t=0.25$. To the left: explosion, to the right: implosion. The numerical solutions (circles) are computed using $N_x=N_y=100$ elements.\label{fig:explosion_implosion}}
    \end{center}
  \end{figure}

\begin{figure}
    \begin{center}
      \psfrag{XXX}[][][0.7][0]{position}
      \psfrag{YYY}[][][0.7][0]{density} 
      \includegraphics[width = 0.45\textwidth]{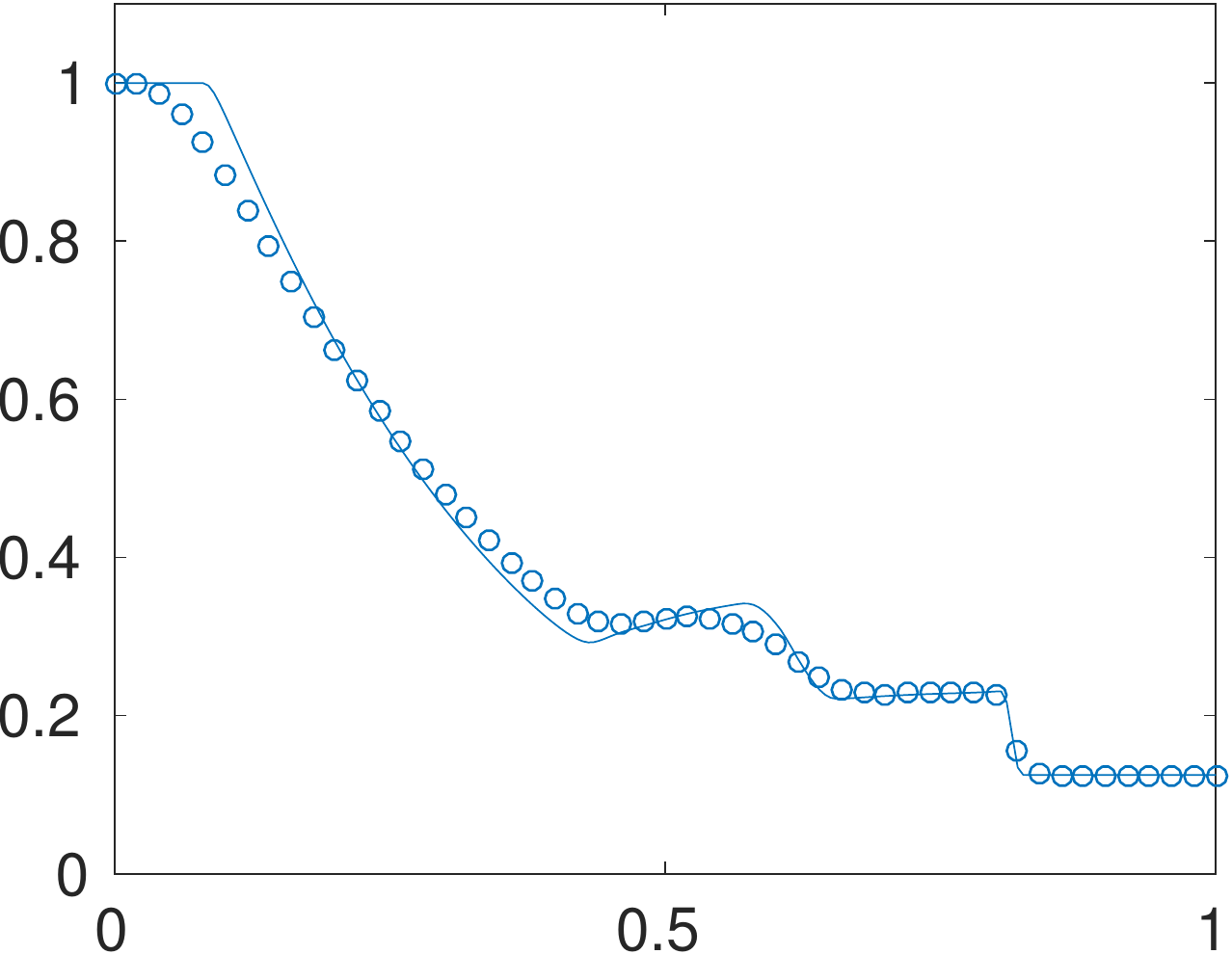}\hspace{0.5cm}
      \includegraphics[width = 0.45\textwidth]{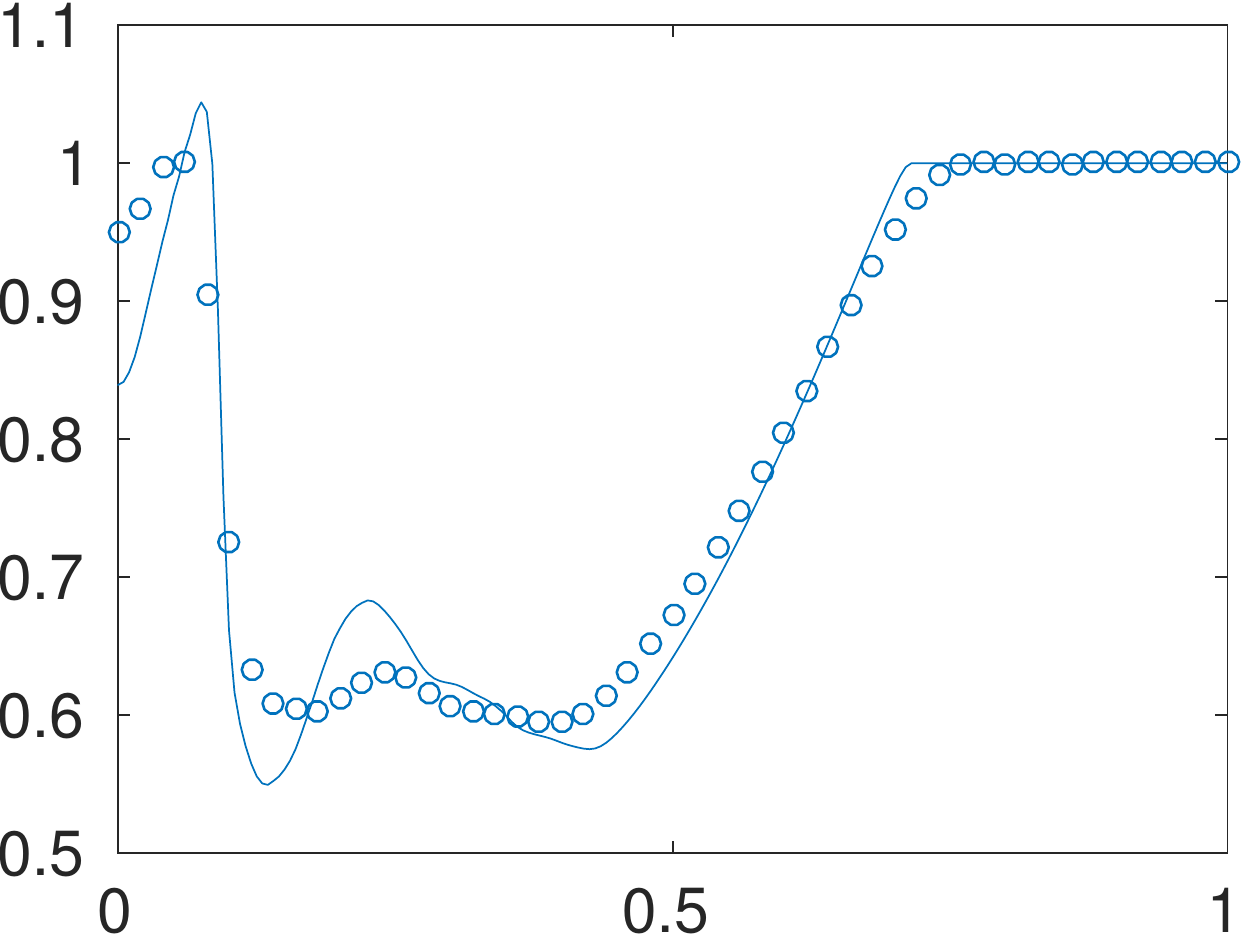}
    \caption{Cross section of the density for radially symmetric problem along $x$ axis at time $t=0.25$. To the left: explosion, to the right: implosion. The numerical solutions (circles) are computed using $N_x=N_y=100$ elements, ``exact'' solutions (solid lines) are computed using $N_x=N_y=400$ elements.\label{fig:radially_symmetric}}
    \end{center}
  \end{figure}

The boundary conditions are imposed by setting the solution on the boundary so that it stays the same as the solution at the initial time. The simulation is performed until time $t=0.25$, before any waves reach the boundary of the domain. We plot the 2D solution in Figure (\ref{fig:explosion_implosion}). In Figure (\ref{fig:radially_symmetric}), we present a cross section of the density at time $t=0.25$ with $N_x=N_y=100$ elements against computed ``exact'' solution obtained with $N_x=N_y=400$ elements.

\subsubsection{Shock Vortex Interaction}
Next we consider the interaction of a shock and a vortex. In general shock-vortex interactions can produce small scales in the form of acoustic waves, and other interesting wave patterns. It has received substantial interest in the literature, see for example \cite{rault2003shock,zhang2005multistage,denet2015model,dumbser2014posteriori}. 

In this experiment, a strong stationary shock with Mach number $2/\sqrt{1.4} \approx 1.69$ collides with a weak vortex with a Mach number $6/2\pi \approx 0.81$. The computational  domain is $D=[-9,3]\times [-4,4]$ and the initial data is  
\begin{equation}
    (\rho,u,v,p)(r,t) =
  \begin{cases}
    (\rho_{vor},u_{vor},v_{vor},p_{vor})             & x>-4, \\
    (2.18,-0.92,0,3.17)                          & x<-4, \\
  \end{cases} 
\end{equation}
where
\begin{align}
  \rho_{vor} & = \left[1-\frac{(\gamma-1)\beta^2}{8\gamma\pi^2}e^{1-x^2-y^2}\right]^{1/(\gamma-1)}\\
u_{vor}    & = 2-\frac{\beta}{2\pi}ye^{(1-x^2-y^2)/2} \\
v_{vor}    & = \frac{\beta}{2\pi}xe^{(1-x^2-y^2)/2} \\
p_{vor}    & = \rho^{\gamma},
\end{align}
and $\beta=6$.

\begin{figure}
  \begin{center}
    \includegraphics[width = 0.45\textwidth]{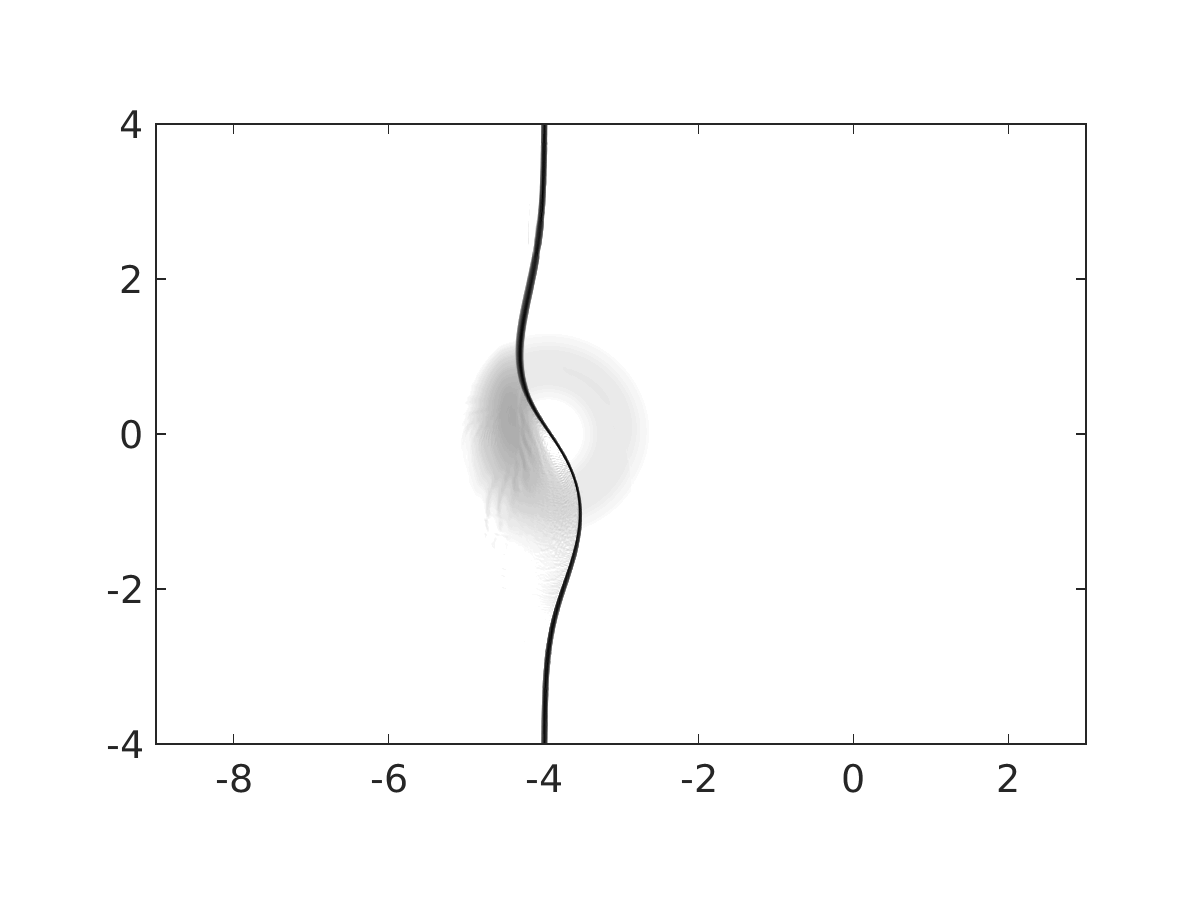}\hspace{0.2cm}
    \includegraphics[width = 0.45\textwidth]{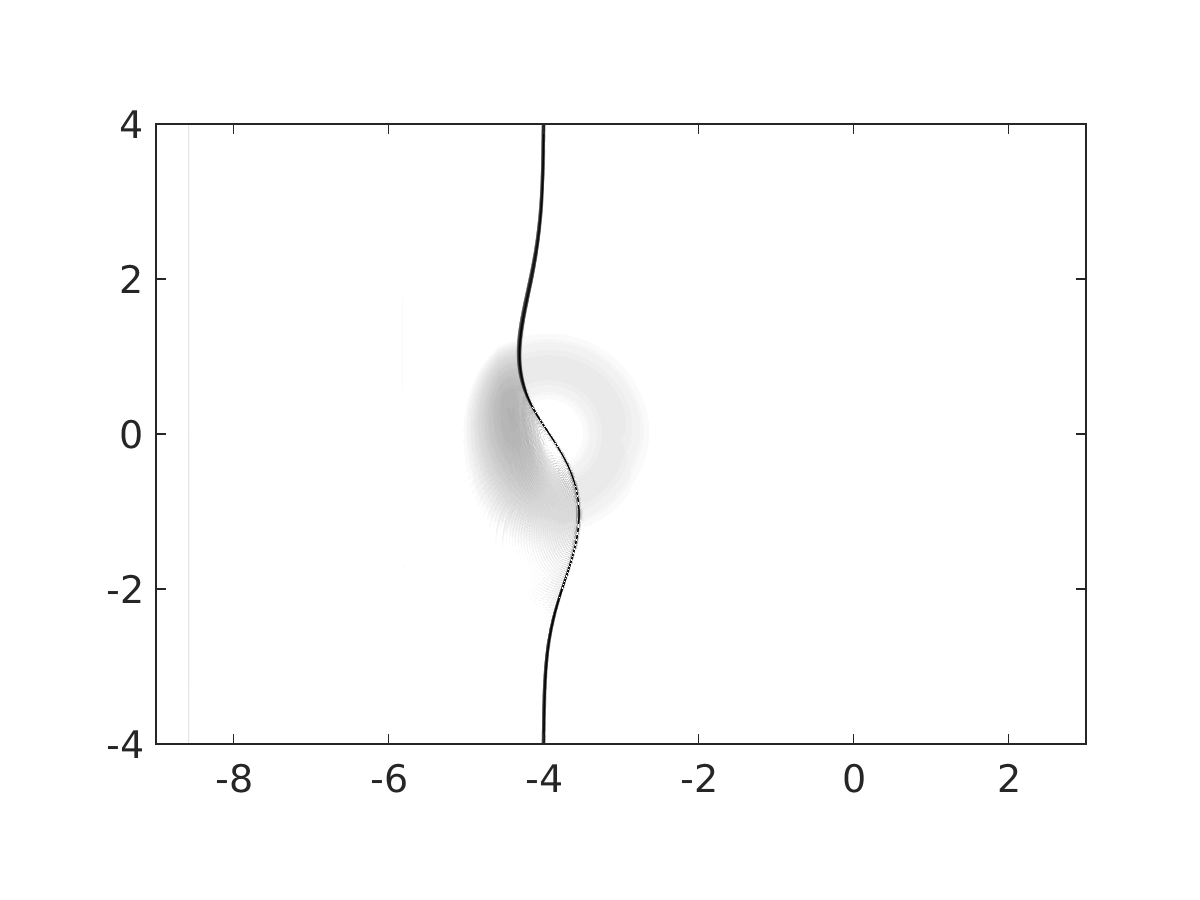}
    \includegraphics[width = 0.45\textwidth]{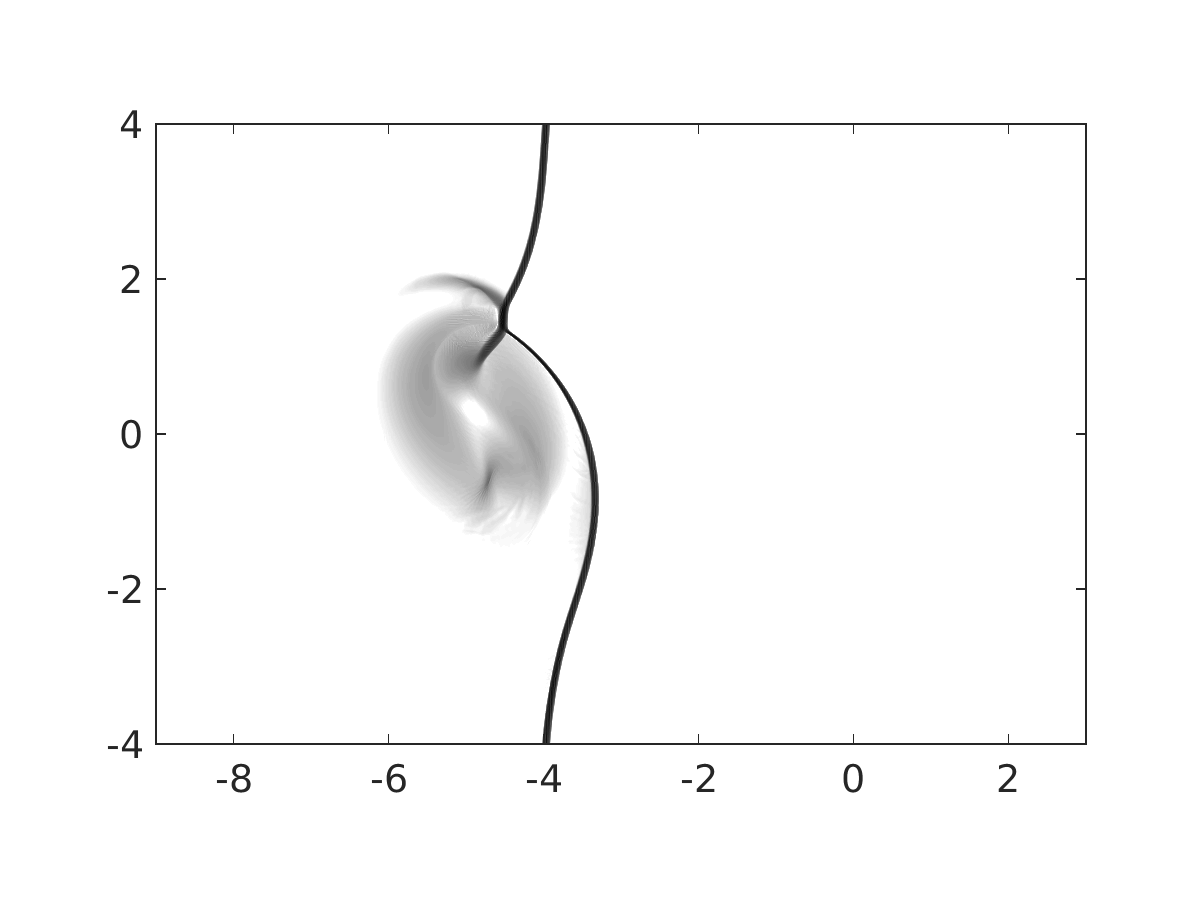}\hspace{0.2cm}
    \includegraphics[width = 0.45\textwidth]{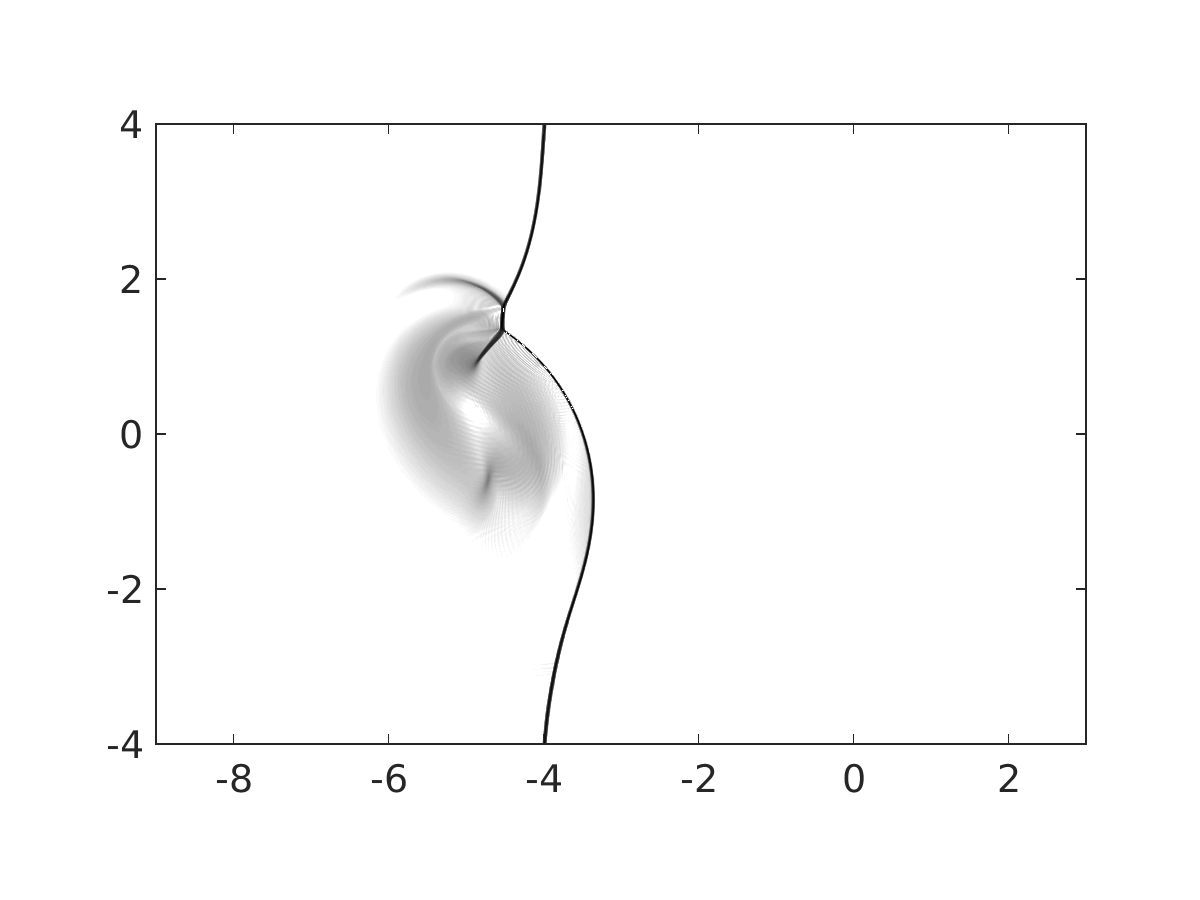}
    \includegraphics[width = 0.45\textwidth]{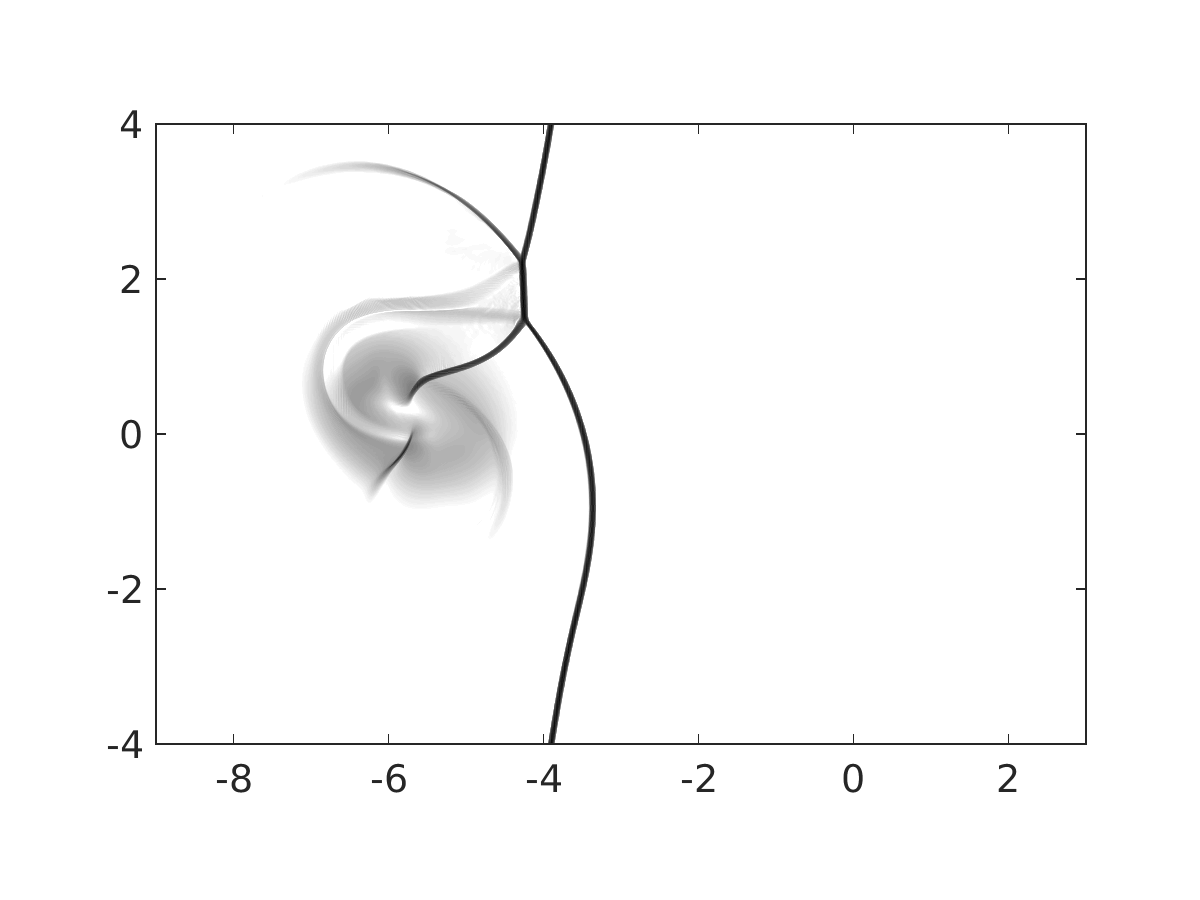}\hspace{0.2cm}
    \includegraphics[width = 0.45\textwidth]{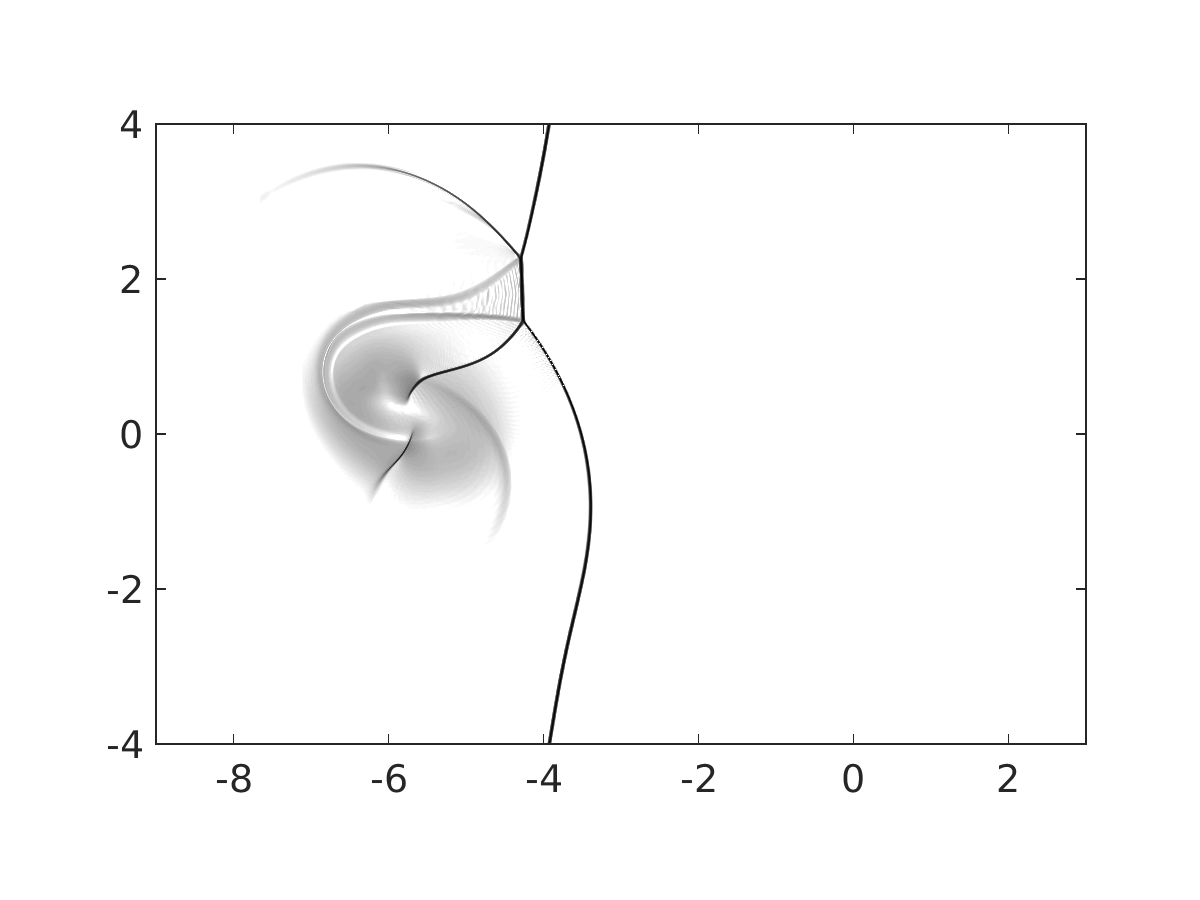}
    \includegraphics[width = 0.45\textwidth]{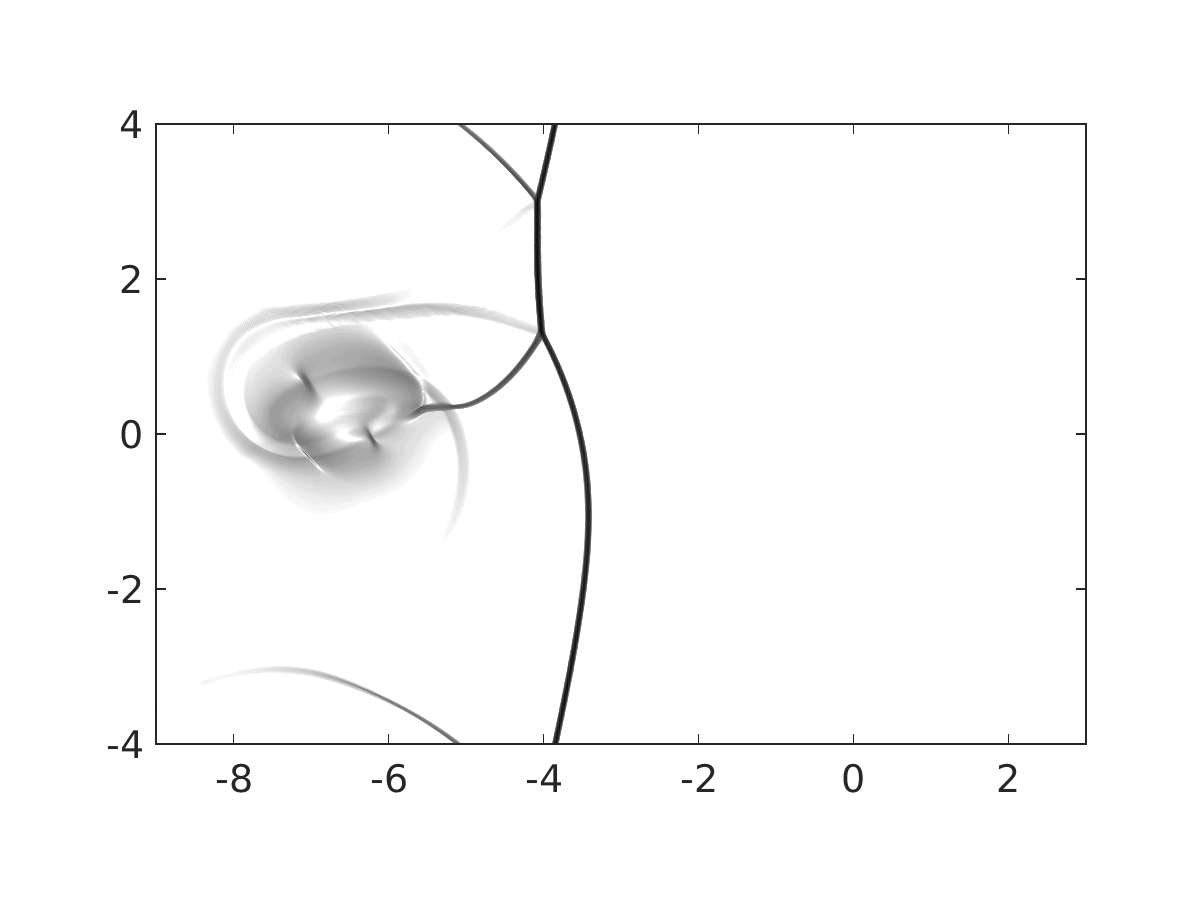}\hspace{0.2cm}
    \includegraphics[width = 0.45\textwidth]{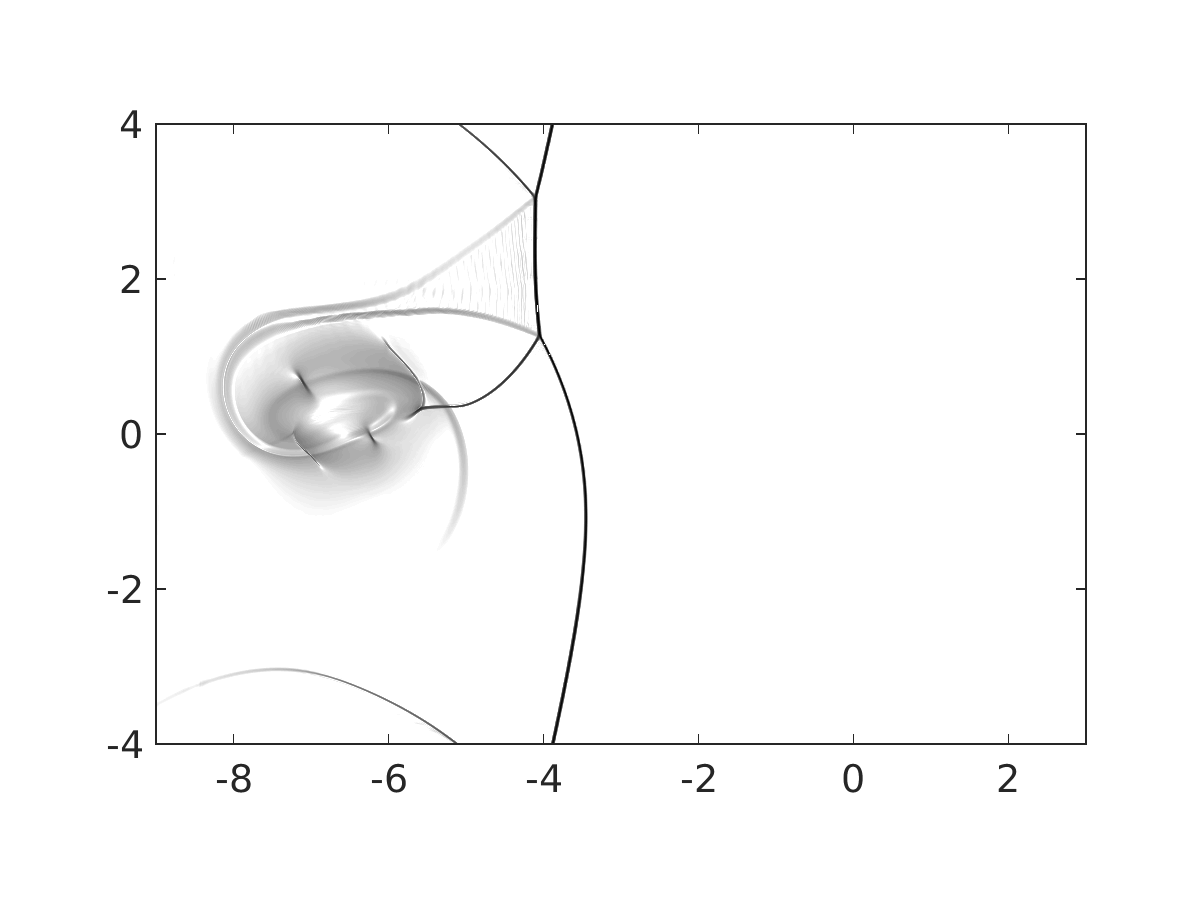}

    \caption{The density schlieren at different times, from top to bottom $t\approx 1.97$, $t\approx 2.95$, $t\approx 3.94$ and $t\approx 4.92$. Left: vortex shock interaction 1, right: vortex interaction 2, with parameters given in Table \ref{tab:parameters_2D}. The numerical solutions are obtained with $N_x=720$, $N_y=480$ elements. \label{fig:schil_vortex1}}
  \end{center}
\end{figure}

As the vortex passes through the shock, the shock is distorted and the vortex is compressed into an elliptical shape. This phenomena is due to the fact that the vortex is relatively weak compared to the shock. The results are consistent with \cite{zhang2005multistage}. In Figure \ref{fig:schil_vortex1}, we compare snapshots of the density Schlieren using two different sets of entropy viscosity parameters, see Table \ref{tab:parameters_2D}. Although the schlierens are plotted on the same color scale, notice that the structures are more pronounced in the pictures on the right column.

\subsubsection{Fluid Flow in Jet}
As a final example we simulate a planar Mach $2/\sqrt{1.4} \approx 1.69$ jet. The domain $D=[-15,55]\times[-17.5,17.5]$ is discretized using of $500 \times 250$ cells. The initial data is 
\begin{equation}
    (\rho,u,v,p)(x,y,t) = (1,0,0,1).      
\end{equation}

We model the jet nozzle by a simple momentum forcing over a $1 \times 1$ patch at the left edge of the computational. The jet is started impulsively causing a relatively strong compression to be generated. This wave sharpens up to a shock wave that is handled by the entropy viscosity as it is propagated from the nozzle and out into a damping absorbing layer of super-grid type, see \cite{AppCol08}. 

\begin{figure}
  \begin{center}
    \includegraphics[width = 0.45\textwidth]{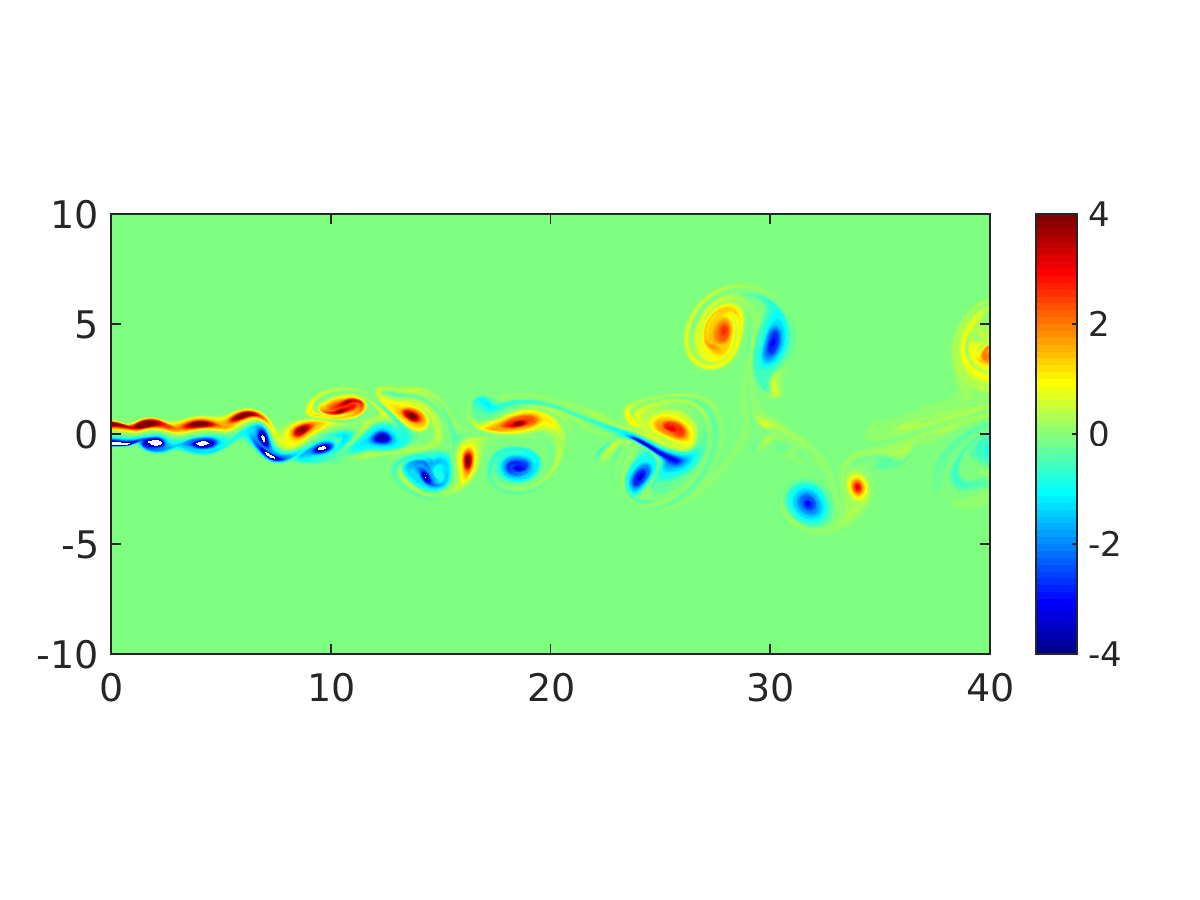}\hspace{0.2cm}
    \includegraphics[width = 0.45\textwidth]{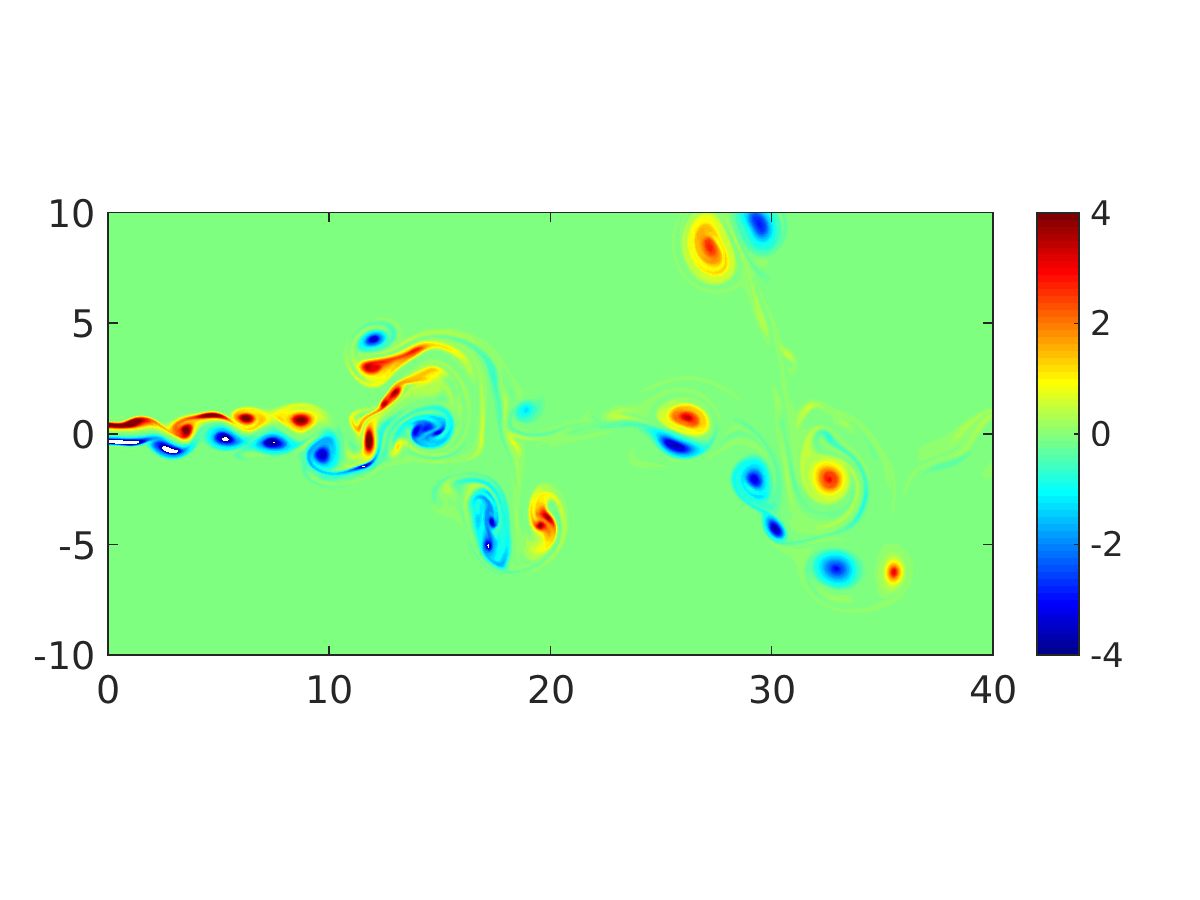}
    \includegraphics[width = 0.45\textwidth]{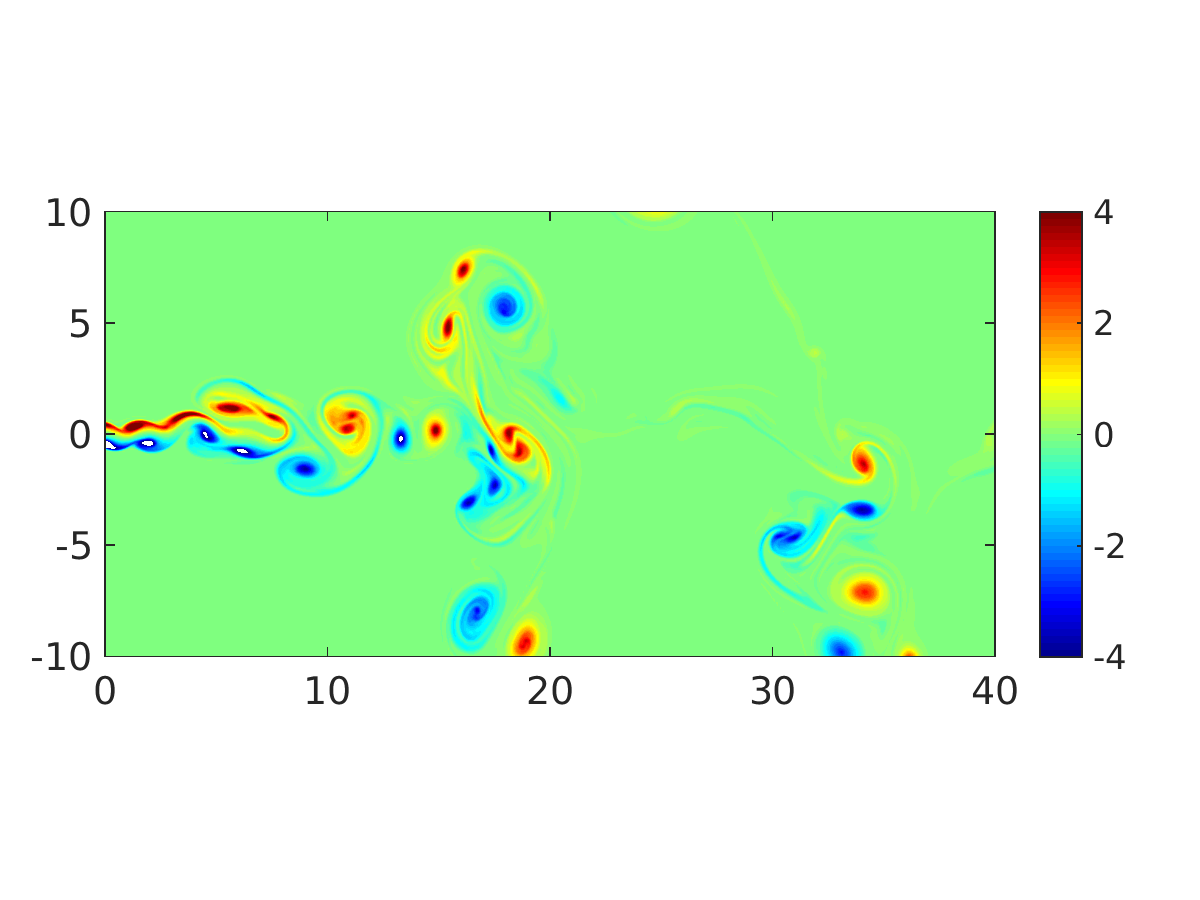}\hspace{0.2cm}
    \includegraphics[width = 0.45\textwidth]{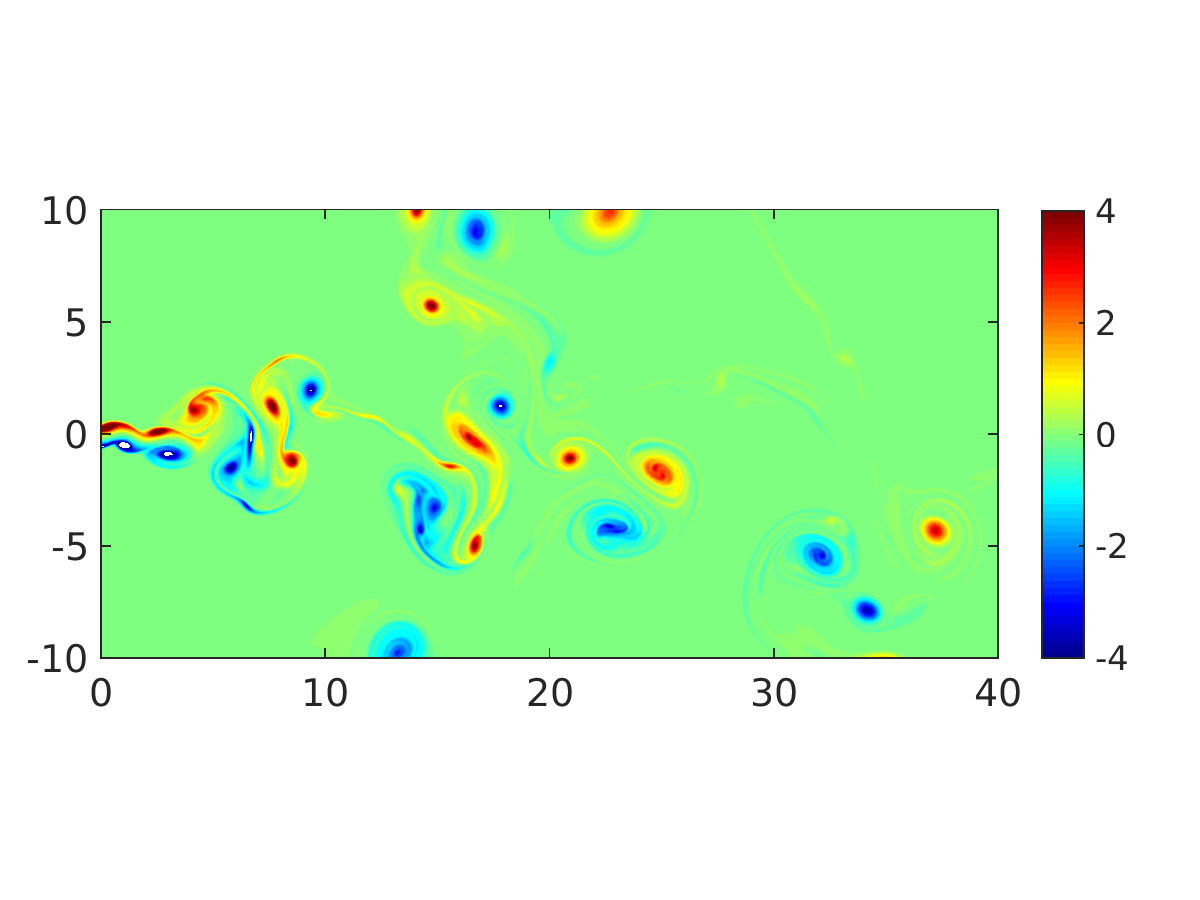}
    \caption{The vorticity at different times, from left to right, top to bottom $t=103.26$, $t=118.93$, $t=134.89$, and $t=150.14$. The numerical solutions are obtained with $N_x=500$, $N_y=250$ elements. \label{fig:vorz_jet}}
    \end{center}
\end{figure}

\begin{figure}
  \begin{center}
    \includegraphics[width = 0.45\textwidth]{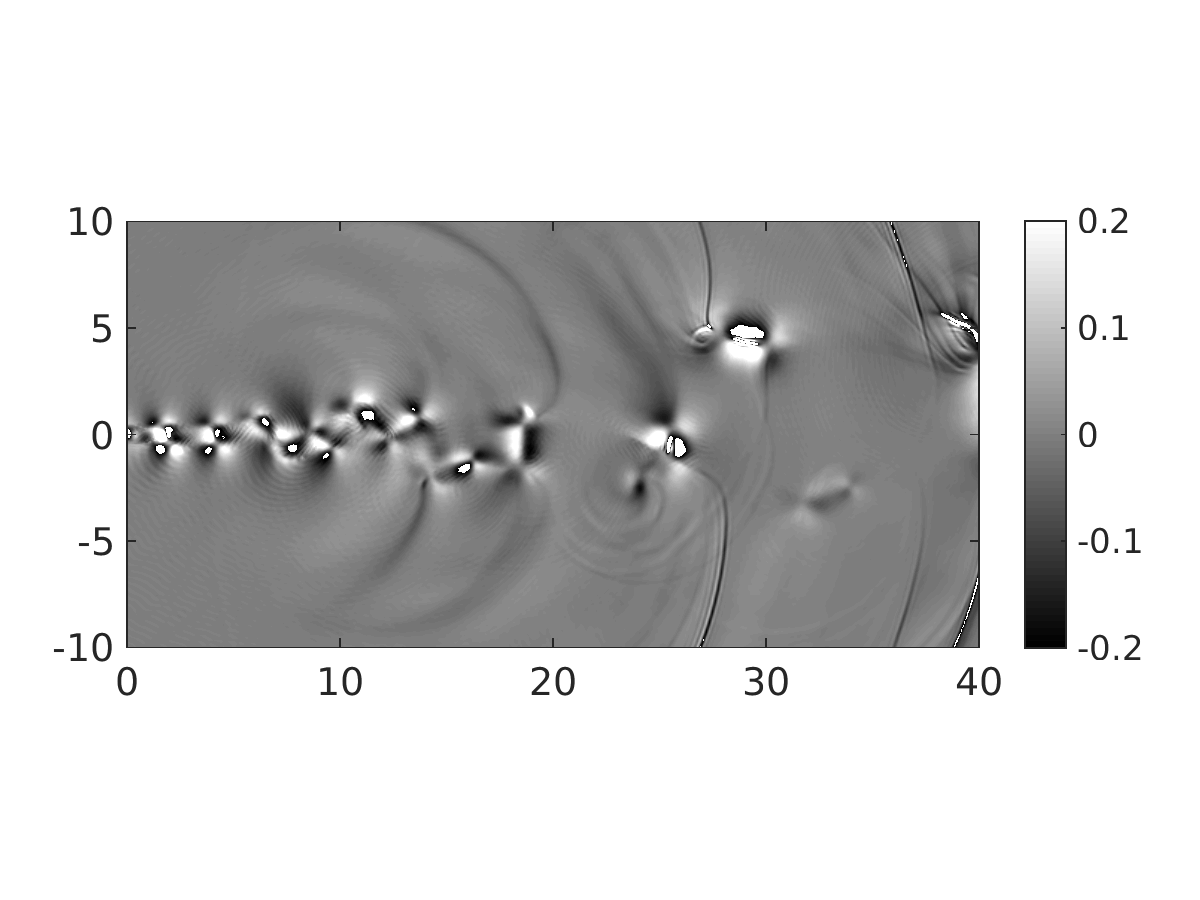}\hspace{0.2cm}
    \includegraphics[width = 0.45\textwidth]{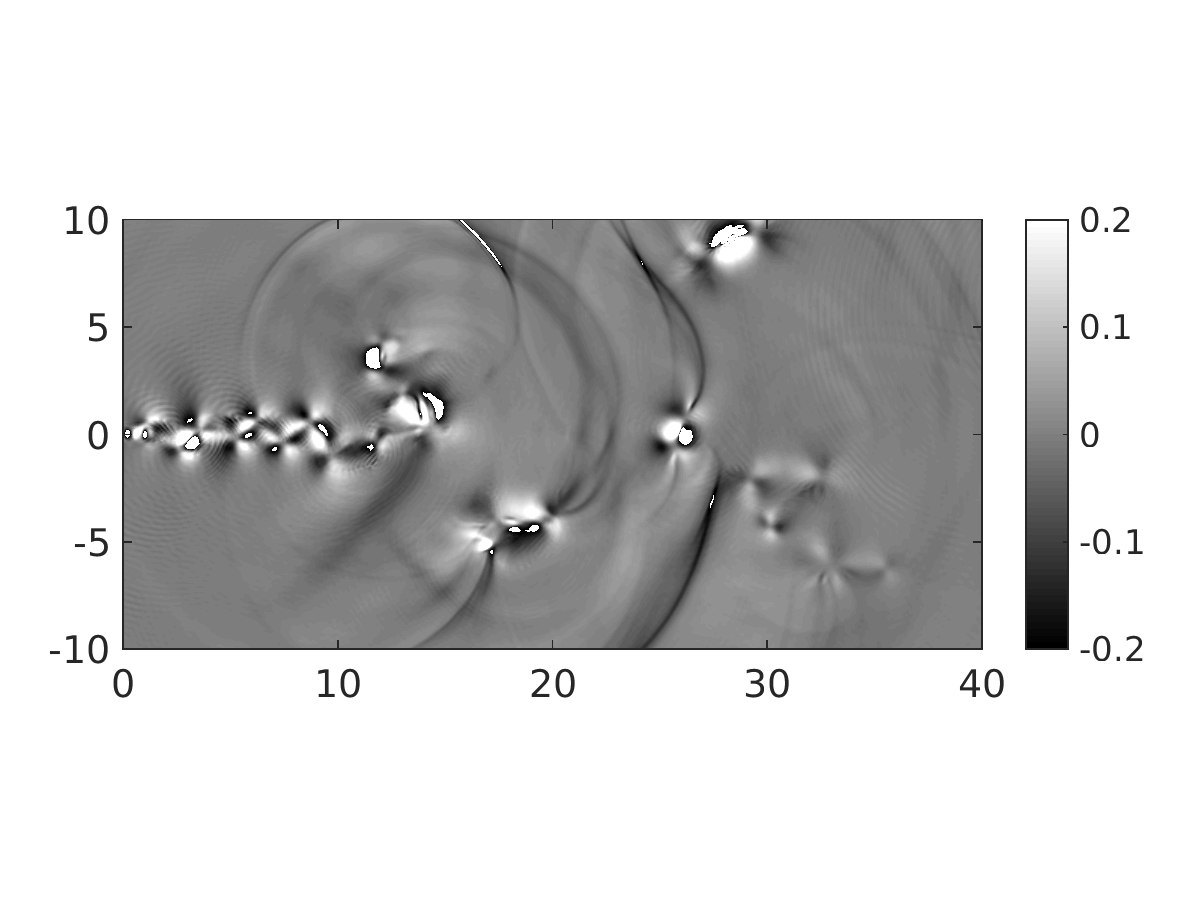}
    \includegraphics[width = 0.45\textwidth]{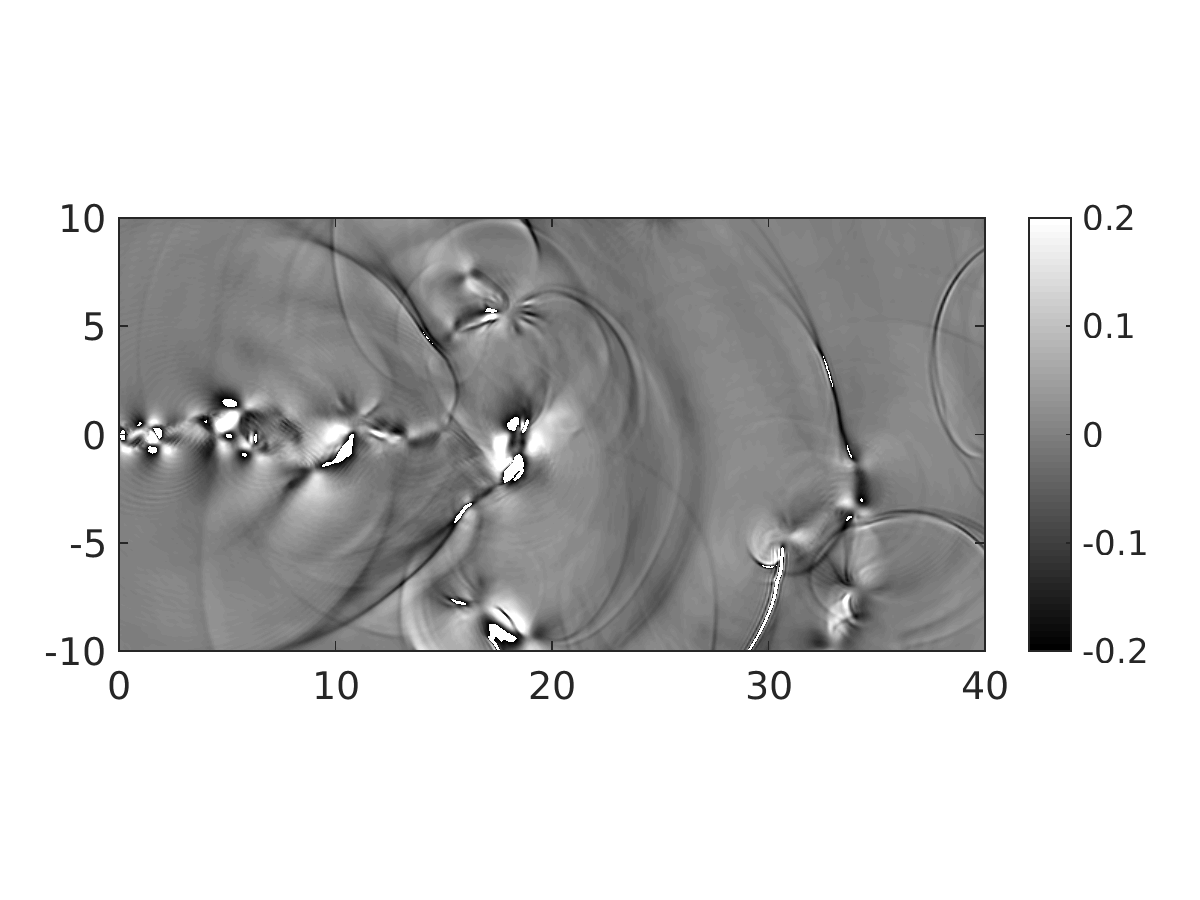}\hspace{0.2cm}
    \includegraphics[width = 0.45\textwidth]{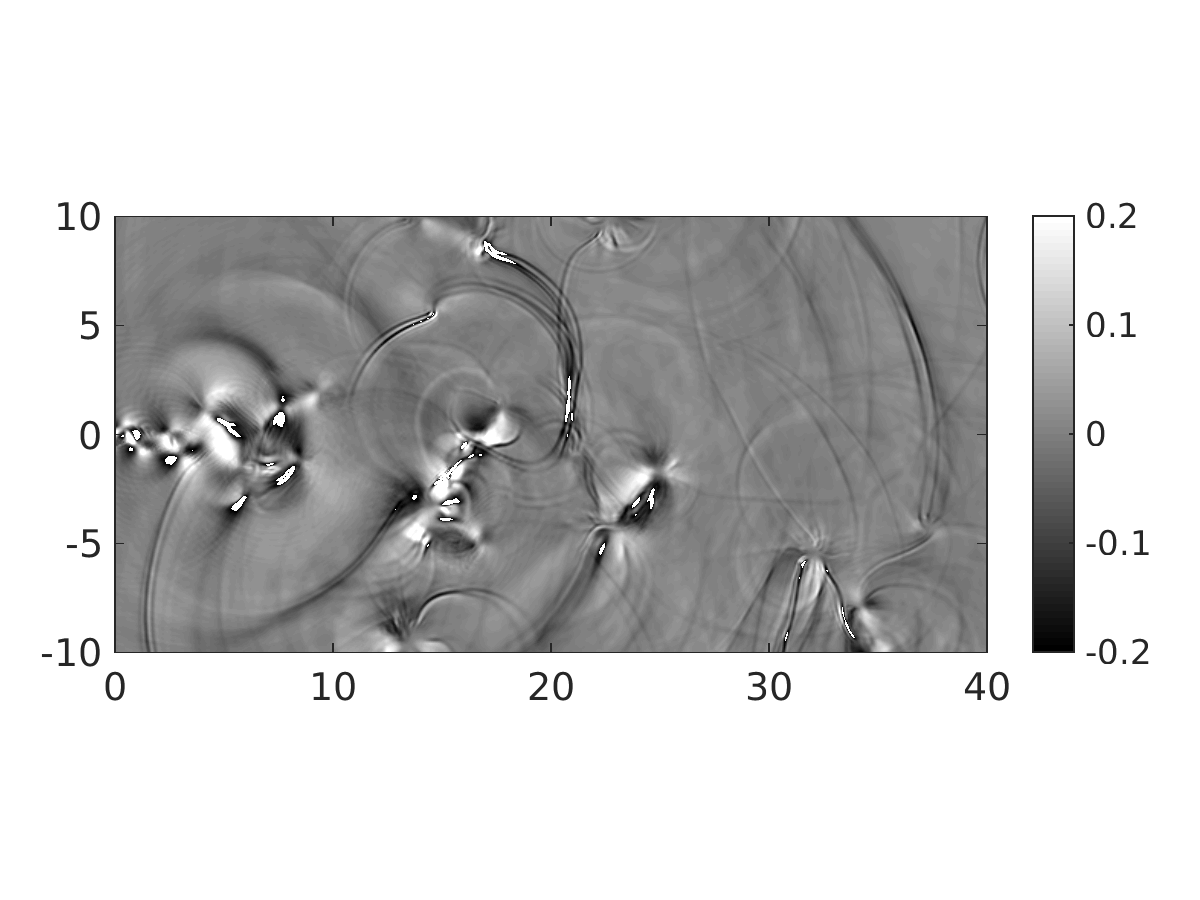}
    \caption{The dilatation at different times, from left to right, top to bottom $t=103.26$, $t=118.93$, $t=134.89$, and $t=150.14$. The numerical solutions are obtained with $N_x=500$, $N_y=250$ elements. \label{fig:dil_jet}}
    \end{center}
\end{figure}

\begin{figure}
  \begin{center}
    \includegraphics[width = 0.45\textwidth]{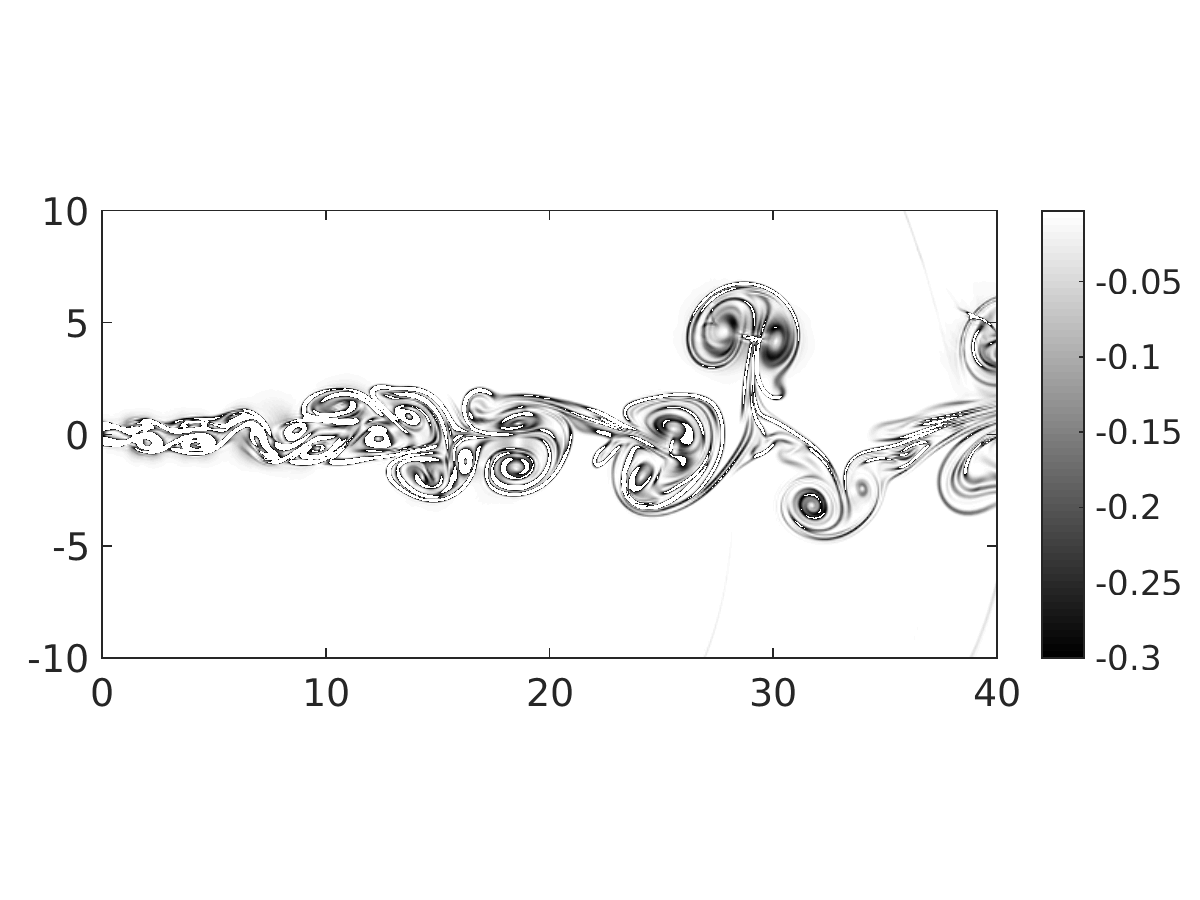}\hspace{0.2cm}
    \includegraphics[width = 0.45\textwidth]{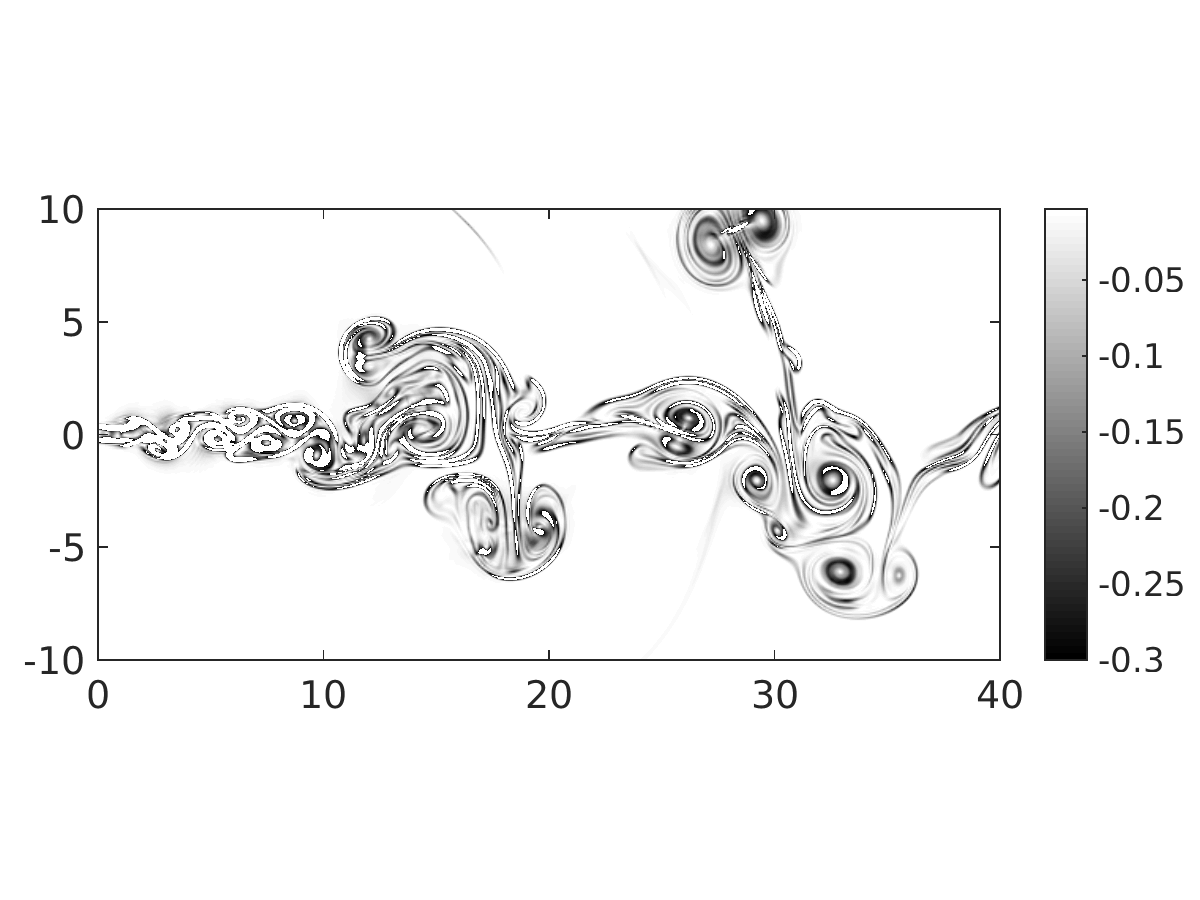}
    \includegraphics[width = 0.45\textwidth]{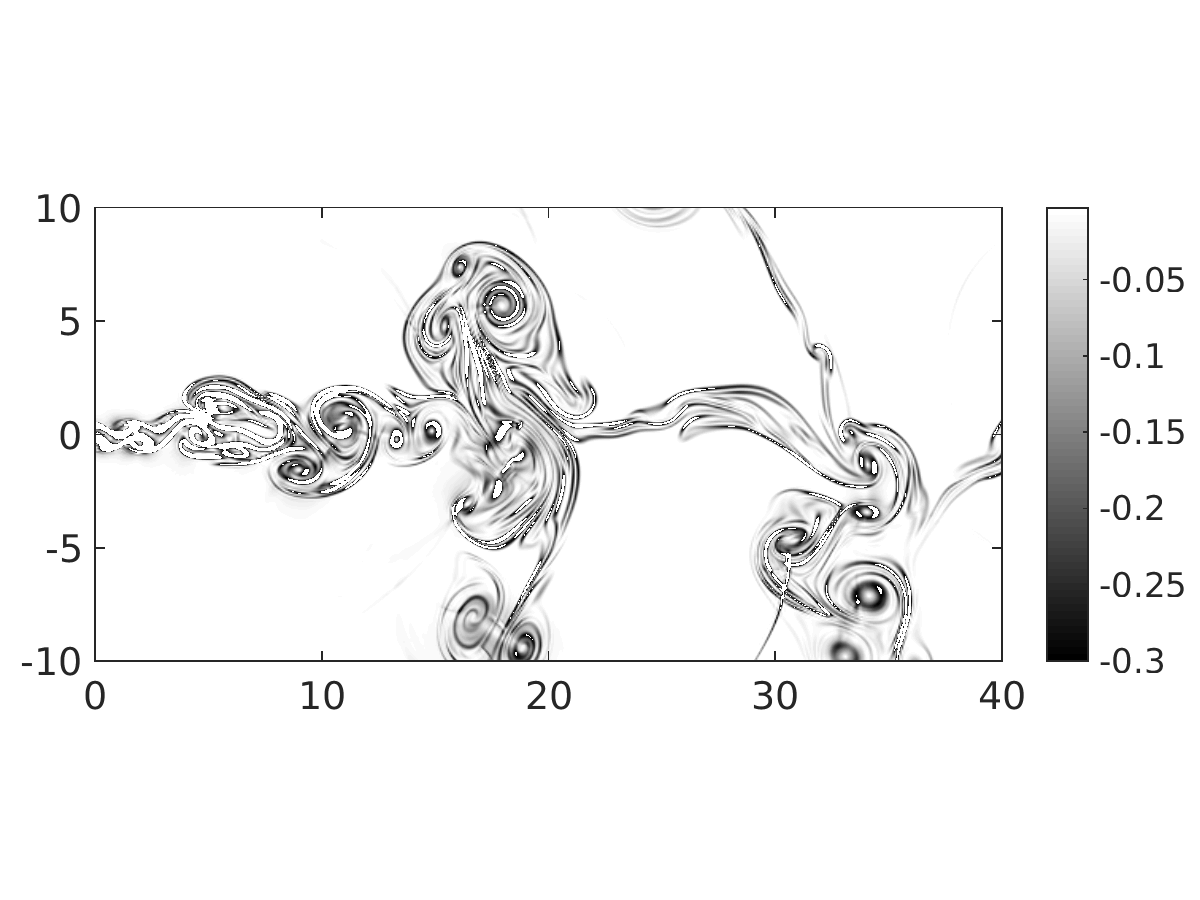}\hspace{0.2cm}
    \includegraphics[width = 0.45\textwidth]{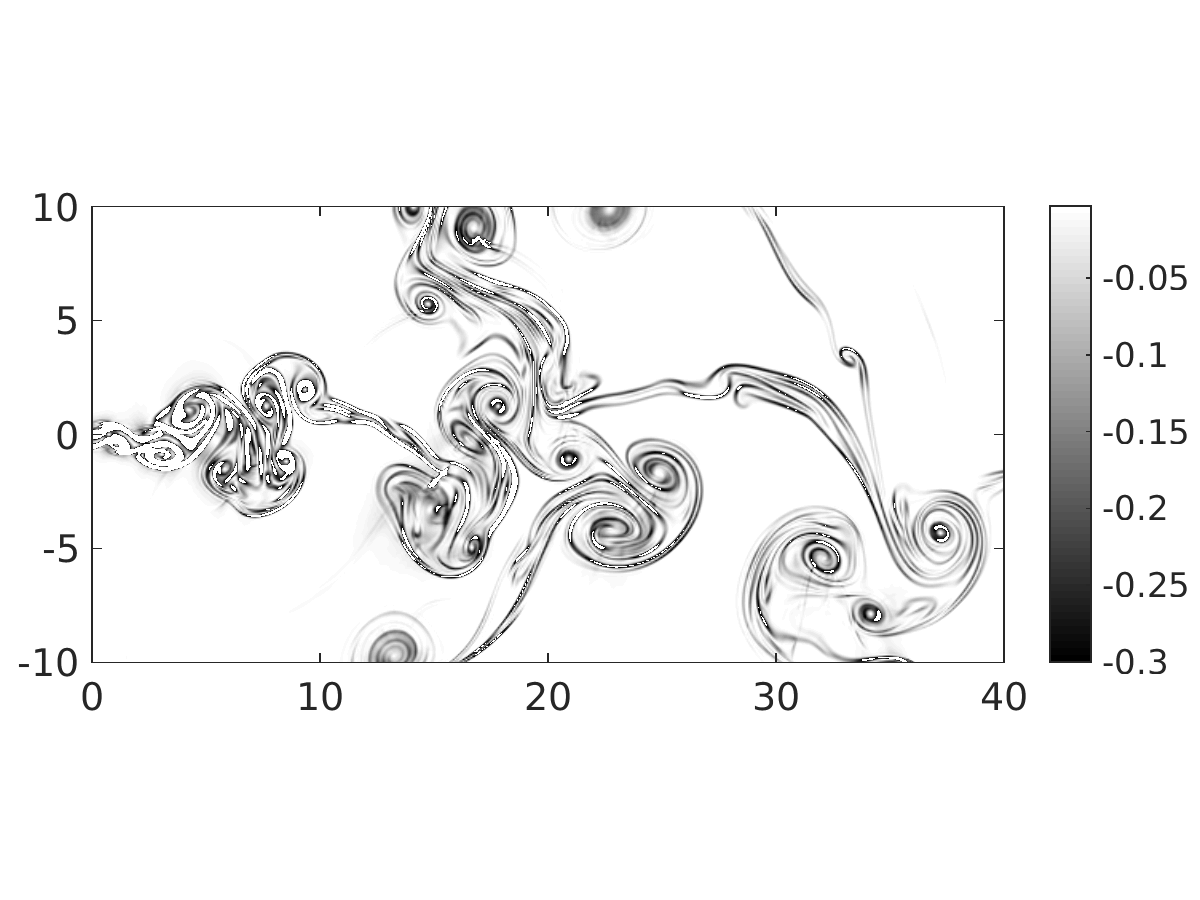}
    \caption{The density schlieren at different times, from left to right, top to bottom $t=103.26$, $t=118.93$, $t=134.89$, and $t=150.14$. The numerical solutions are obtained with $N_x=500$, $N_y=250$ elements. \label{fig:schil_jet}}
  \end{center}
\end{figure}

In Figures (\ref{fig:vorz_jet})-(\ref{fig:schil_jet}) we display snapshots of the vorticity, dilatation and density Schlieren. We note that the viscosity we use here is purely for the regularization of shocks so there is no reason to believe that the flow that we compute resembles reality. Nevertheless, the example illustrates the new methods ability to handle rapidly started flows. Also, it is likely that the particular form of the artificial viscosity does not effect the robustness of the method. 

  \section{Conclusion}
In conclusion we have demonstrated that flux-conservative Hermite methods are suitable for solving nonlinear conservation laws, especially in the presence of shocks. The new methods still converges at a rate of $(2m+1)$ for smooth problems.

The adaptation of the entropy viscosity method to Hermite methods successfully suppresses oscillations near shocks, but  we find that our current implementation is quite dissipative when solving the Shu-Osher problem. For contact waves we proposed a modification to the entropy viscosity method which eliminates a large amount of the  spurious viscosity at contact discontinuities.

\bibliographystyle{plain}
\bibliography{kornelus}

\end{document}